\documentclass{amsart}
\newtheorem{thm}{Theorem}
\newtheorem{cor}[thm]{Corollary}
\newtheorem{lem}[thm]{Lemma}

\newtheorem{prop}[thm]{Proposition}

\newtheorem{rem}[thm]{Remark}

\newtheorem{defn}[thm]{Definition}
\newtheorem{notation}[thm]{Notation}
\numberwithin{equation}{section}

\def \R{ \bf R \rm }
\def \p{ \pi }
\def \r{ \rho }
\def \t{ \theta }
\def \l{ \lambda }
\def \e{ \varepsilon }
\def \f{ \varphi }
\def \a{ \alpha }
\def \s{ \sigma }
\def \j{ \chi }
\def \m{ \mu }
\def \N{ \bf N \rm }
\def \d{ \delta}
\def \Pf{ \partial _Pf }
\def \F{ \Phi }
\def \z{ \zeta }
\def \g{ \gamma }
\begin{document}

\title[Proximal calculus on Riemannian manifolds]{Proximal calculus on
Riemannian manifolds, with applications to fixed point theory}
\author{Daniel Azagra and Juan Ferrera}

\address{Departamento de An{\'a}lisis Matem{\'a}tico\\ Facultad de
Matem{\'a}ticas\\ Universidad Complutense\\ 28040 Madrid, Spain}

\date{March 24, 2004}

\email{azagra@mat.ucm.es, ferrera@mat.ucm.es}



\begin{abstract}
We introduce a proximal subdifferential and develop a calculus
for nonsmooth functions defined on any Riemannian manifold $M$.
We give several applications of this theory, concerning:
\begin{enumerate}
\item differentiability and geometrical
properties of the distance function to a closed subset $C$ of
$M$;
\item solvability and implicit function theorems for
nonsmooth functions on $M$;
\item conditions on the existence of a circumcenter for three
different points of $M$; and especially
\item fixed point theorems for expansive and nonexpansive mappings and
certain perturbations of such mappings defined on $M$.
\end{enumerate}
\end{abstract}

\maketitle

\section{Introduction and main results}

The proximal subdifferential of lower semicontinuous real-valued
functions is a very powerful tool that has been extensively
studied and used in problems of optimization, control theory,
differential inclusions, Lyapounov Theory, stabilization, and
Hamilton-Jacobi equations; see \cite{C} and the references
therein.

In this paper we will introduce a notion of proximal
subdifferential for functions defined on a Riemannian manifold $M$
(either finite or infinite dimensional) and we will develop the
rudiments of a calculus for nonsmooth functions defined on $M$.
Next we will prove an important result concerning inf-convolutions
of lower semicontinuous functions with squared distance functions
on $M$, from which a number of interesting consequences are
deduced. For instance, we show a Borwein-Preiss type variational
principle for lower semicontinuous functions defined on $M$, and
we study some differentiability and geometrical properties of the
distance function to a closed subset $C$ of $M$. Then we establish
a Decrease Principle from which we will deduce some Solvability
and Implicit Function Theorems for nonsmooth functions on $M$. We
give two applications of the Solvability Theorem. First, we give
some conditions guaranteeing the existence of a circumcenter for
three given points of $M$. Second and most important, we provide
several fixed point theorems for expansive and nonexpansive
mappings and certain perturbations of such mappings defined on
$M$. Observe that, in general, very small perturbations of
mappings having fixed points may lose them: consider for instance
$f:\mathbb{R}\to\mathbb{R}$, $f(x)=x+\varepsilon$. Note also that
most of the known fixed point theorems (such as Brouwer's, Lefschetz's,
Schauder's or the Banach contraction mapping principle) rely either on compactness
or on contractiveness, see \cite{Brown, Cronin} for instance.
However, our results hold for (possibly expansive) mappings on
(possibly noncompact) complete Riemannian manifolds.

Let us give a brief sample of the corollaries on fixed points that
we will be deducing from our main theorems in the last section of
this paper.
%%%%%%%%%%%%%%%%%
\begin{cor}
Let $M$ be a complete Riemannian manifold with a positive
injectivity radius $\rho=i(M)$.  Let $x_{0}$ be a fixed point of
a $C^1$ smooth mapping $G:M\to M$ such that $G$ is $C$-Lipschitz
on a ball $B(x_{0},R)$. Let $H:M\to M$ be a differentiable
mapping. Assume that $0<2R\max\{1, C\}<\rho$, that $$\langle
L_{xH(G((x)))}h, L_{G(x)H(G(x))}dG(x)(h)\rangle_{F(x)}\leq K<1$$
for all $x\in B(x_{0},R)$ and $h\in TM_{x}$ with $\|h\|_{x}=1$,
and that $\|dH(y)-L_{yH(y)}\|<\varepsilon/C$ for every $y\in
G(B(x_{0}), R))$, where $\varepsilon<1-K$, and $d(x_{0},
H(G(x_{0})))< R(1-K-\varepsilon)$. Then $F=H\circ G$ has a fixed
point in $B(x_{0}, R)$.
\end{cor}
%%%%%%%%%%%%%%%%%
Here $L_{xy}$ stands for the parallel transport along the (unique
in this setting) minimizing geodesic joining the points $x$ and
$y$. The hypotheses on $H$ mean that $H$ is relatively close to
the identity, so the perturbation brought on $G$ by its
composition with $H$ is relatively small.
%%%%%%%%%%%%%%%%%
\begin{cor}
Let $(M,+)$ be a complete Riemannian manifold with a Lie group
structure. Let $x_0$ be a fixed point of a differentiable
function $G:M\to M$ satisfying the following condition: $$\langle
h, dG(x_0)(h)\rangle \leq K<1 \ \ \ \hbox{for every} \ \ \
\|h\|=1.$$ Then there exists a positive $\delta$ such that for
every Lipschitz mapping $H:M\to M$ with Lipschitz constant
smaller than $\delta$, the mapping $G+H:M\to M$ has a fixed point.
\end{cor}
%%%%%%%%%%%%%%%%%
The condition on the differential of $G$ is satisfied, for
instance, if $G$ locally behaves like a multiple of a rotation
round the point $x_{0}$, but notice that $G$ may well be
expansive. Consider for instance
$G:\mathbb{R}^{2}\to\mathbb{R}^{2}$, $G(x,y)=23(y,-x)$; in this
case we can take $K=0$, but $G$ is clearly expansive.

It should be stressed that these fixed point results are new even
in the case when $M=\mathbb{R}^{n}$ or any Hilbert space. In
particular we have the following.
%%%%%%%%%%%%%%%%%
\begin{cor}
Let $X$ be a Hilbert space, and let $x_0$ be a fixed point of a
differentiable mapping $G:X\to X$ satisfying the following
condition: $$\langle h,DG(x)(h)\rangle\leq K<1 \ \ \ \hbox{for
every} \ \ \ x\in B(x_0,R)\ \ \hbox{and} \ \ ||h||=1.$$ Then we
have that:
\begin{enumerate}
\item If $H$ is a differentiable $L$-Lipschitz mapping, with
$L<1-K$, then $G+H$ has a fixed point in $B(x_0, R)$, provided that
$\|H(x_0)\|<R(1-K-L)$.
\item If $H:X\to X$ is a differentiable mapping such that
$\|DH(G(x))-I\|<\varepsilon$ for every $x\in B(x_{0}, R)$, then
$F=H\circ G$ has a fixed point in $B(x_{0}, R)$, provided that
$K+\varepsilon<1$ and $\|H(x_{0})-x_{0}\|<R(1-K-\varepsilon)$.
\end{enumerate}
\end{cor}
%%%%%%%%%%%%%%%%%
All of these results and many other things will be proved in
section 3.

\medskip

This paper should be compared with \cite{AFL2}, where a theory of
viscosity subdifferentials for functions defined on Riemaniann
manifolds is established and applied to show existence and
uniqueness of viscosity solutions to Hamilton-Jacobi equations on
such manifolds.

Let us recall the definition of the proximal subdifferential for
functions defined on a Hilbert space $X$. A vector $\zeta{\in X}$
is called a {\sl proximal subgradient} of a lower semicontinuous
function $f$ at $x\in dom f:=\{y\in X: f(y)<+\infty\}$ provided
there exist positive numbers $\sigma$ and $\eta$ such that
\[
f\left( y \right) \geq f\left( x \right) + \left\langle {\zeta ,y
- x} \right\rangle  - \sigma \left\| {y - x} \right\|^2\textrm{
for all } y \in B\left( {x,\eta } \right).
\]
The set of all such $\zeta$ is denoted $\partial_{P}f(x)$, and is
referred to as the {\sl proximal sub\-differential}, or
P-subdiferential. A comprehensive study of this subdifferential
and its numerous applications can be found in \cite{C}.

Before giving the definition of proximal subdifferential for a
function defined on a Riemannian manifold, we must establish a few
preliminary results.

The following result is proved in \cite[Corollary 2.4]{AFL1}
%%%%%%%%%%%%%%%%%
\begin{prop}\label{characterization of the proximal subdiferential through C2 functions}
Let $X$ be a real Hilbert space, and $f:X\longrightarrow
(-\infty,\infty]$ be a proper, lower semicontinuous function.
Then,
    $$
    \partial_{p}f(x)=\{\varphi'(x) : \varphi\in
    C^{2}(X,\mathbb{R}),\,
    f-\varphi \textrm{ attains a local minimum at } x\}.
    $$
\end{prop}
%%%%%%%%%%%%%%%%%
In particular this implies that $\partial_{P}f(x)\subseteq
D^{-}f(x)$, where $D^{-}f(x)$ is the viscosity subdifferential of
$f$ at $x$.
%%%%%%%%%%%%%%%%%
\begin{lem}\label{equivalence...}
Let $X_1$\ and $X_2$\ be two real Hilbert spaces, $\Phi :X_2 \to
X_1$\ a $C^2$ diffeomorphism, $f:X_1 \to (-\infty ,+\infty ]$\ a
lower semicontinuous function. Then $v\in \Pf (x_1)$\ if and only
if $D\F (x_2)^*(v)\in \partial (f\circ \F )(x_2)$, where $\F
(x_2)=x_1$.
\end{lem}
%%%%%%%%%%%%%%%%%
\begin{proof}
This is a trivial consequence of Proposition \ref{characterization
of the proximal subdiferential through C2 functions}, bearing in
mind that compositions with diffeomorphisms preserve local minima.
\end{proof}
%%%%%%%%%%%%%%%%%

%%%%%%%%%%%%%%%%%
\begin{cor}
Let $M$\ be a Riemannian manifold, $p\in M$, $(\f _i, U_i)$\
$i=1,2$, two charts with $p\in U_1\cap U_2$, and $\f _i(p)=x_i$.
Then $\partial _P(f\circ \f _1^{-1})(x_1)\not= \emptyset$\ if and
only if $\partial _P(f\circ \f _2^{-1})(x_2)\not= \emptyset$.
Moreover, $D(\varphi_{1}\circ\varphi_{2}^{-1})(x_{2})^{*}
(\partial_{P}(f\circ\varphi_{1}^{-1})(x_{1}))=
\partial_{P}(f\circ\varphi_{2}^{-1})(x_{2})$.
\end{cor}
%%%%%%%%%%%%%%%%%

Now we can extend the notion of P-subdifferential to functions
defined on a Riemannian manifold.
%%%%%%%%%%%%%%%%%
\begin{notation}
{\em In the sequel, $M$ will stand for a Riemannian manifold
defined on a real Hilbert space $X$ (either finite dimensional or
infinite-dimensional). As usual, for a point $p\in M$, $TM_{p}$
will denote the tangent space of $M$ at $p$, and
$\exp_{p}:TM_{p}\to M$ will stand for the exponential function at
$p$.

We will also make extensive use of the parallel transport of
vectors along geodesics. Recall that, for a given curve $\gamma:
I\to M$, numbers  $t_{0}, t_{1}\in I$, and a vector $V_{0}\in
TM_{\gamma(t_{0})}$, there exists a unique parallel vector field
$V(t)$ along $\gamma(t)$ such that $V(t_{0})=V_{0}$. Moreover,
the mapping defined by $V_{0}\mapsto V(t_{1})$ is a linear
isometry between the tangent spaces $TM_{\gamma(t_{0})}$ and
$TM_{\gamma(t_{1})}$, for each $t_{1}\in I$. In the case when
$\gamma$ is a minimizing geodesic and $\gamma(t_{0})=x$,
$\gamma(t_{1})=y$, we will denote this mapping by $L_{xy}$, and
we call it the parallel transport from $TM_{x}$ to $TM_{y}$ along
the curve $\gamma$. See \cite{Klingenberg} for general reference
on these topics.

The parallel transport allows us to measure the length of the
``difference" between vectors (or forms) which are in different
tangent spaces (or in duals of tangent spaces, that is, fibers of
the cotangent bundle), and do so in a natural way. Indeed, let
$\gamma$ be a minimizing geodesic connecting two points $x, y\in
M$, say $\gamma(t_{0})=x, \gamma(t_{1})=y$. Take vectors $v\in
TM_{x}$, $w\in TM_{y}$. Then we can define the distance between
$v$ and $w$ as the number
    $$
    \|v-L_{yx}(w)\|_{x}=
    \|w-L_{xy}(v)\|_{y}
    $$
(this equality holds because $L_{xy}$ is a linear isometry between
the two tangent spaces, with inverse $L_{yx}$). Since the spaces
$T^{*}M_{x}$ and $TM_{x}$ are isometrically identified by the
formula $v=\langle v, \cdot\rangle$, we can obviously use the
same method to measure distances between forms $\zeta\in
T^{*}M_{x}$ and $\eta\in T^{*}M_{y}$ lying on different fibers of
the cotangent bundle.}
\end{notation}
%%%%%%%%%%%%%%%%%

%%%%%%%%%%%%%%%%%
\begin{defn}
Let $M$\ be a Riemannian manifold, $p\in M$, $f:M \to (-\infty
,+\infty ]$\ a lower semicontinuous function. We define the
proximal subdifferential of $f$\ at $p$, denoted by $\Pf
(p)\subset TM_{p}$, as $\partial _P(f\circ exp_p)(0)$ (understood
that $\Pf (p)=\emptyset$ for all $p\notin\textrm{dom}f$).
\end{defn}
%%%%%%%%%%%%%%%%%

The following result is an immediate consequence of Lemma
\ref{equivalence...}.

%%%%%%%%%%%%%%%%%
\begin{prop}
Let $M$\ be a Riemannian manifold, $p\in M$, $(\f , U)$\ a chart,
with $p\in U$, and $f:M \to (-\infty ,+\infty ]$\ a lower
semicontinuous function. Then $$\Pf (p)=D\f (p)^*[\partial
_P(f\circ \f ^{-1})(\f (p)].$$
\end{prop}
%%%%%%%%%%%%%%%%%
As a consequence of the definition of
$\partial_{P}(f\circ\exp_{p})(0)$ we get the following.
%%%%%%%%%%%%%%%%%
\begin{cor}
Let $M$\ be a Riemannian manifold, $p\in M$, $f:M \to (-\infty
,+\infty ]$\ a lower semicontinuous function. Then $\z \in \Pf
(p)$\ if and only if there is a $\s >0$\ such that $$ f(q)\geq
f(p)+\langle\z , exp_p^{-1}(q)\rangle-\s d(p,q)^2$$ for every $q$\
in a neighborhood of $p$.
\end{cor}
%%%%%%%%%%%%%%%%%

We can also define the proximal superdifferential of a function
$f$ from a Hilbert space $X$ into $[-\infty, \infty)$ as follows.
A vector $\zeta{\in X}$ is called a {\sl proximal supergradient}
of an upper semicontinuous function $f$ at $x\in dom f$ if there
are positive numbers $\sigma$ and $\eta$ such that
\[
f\left( y \right) \leq f\left( x \right) + \left\langle {\zeta ,y
- x} \right\rangle  + \sigma \left\| {y - x} \right\|^2\textrm{
for all } y \in B\left( {x,\eta } \right).
\]
and we denote the set of all such $\zeta$ by $\partial^{P}f(x)$,
which we call P-subddiferential of $f$ at $x$.

Now, if $M$ is a Riemannian manifold, $p\in M$, $f:M \to [-\infty
,+\infty)$\ an upper semicontinuous function. We define the
proximal superdifferential of $f$\ at $p$, denoted by
$\partial^{P}f (p)\subset TM_{p}$, as $\partial^{P}(f\circ
exp_p)(0)$. As before, we have that $\z \in \partial^{P}f (p)$ if
and only if there is a $\s >0$\ such that $$ f(q)\leq
f(p)+\langle\z , exp_p^{-1}(q)\rangle+\s d(p,q)^2$$ for every $q$\
in a neighborhood of $p$. It is also clear that $\partial^{P}f
(p)=-\partial_{P}(-f) (p)$.

Before going into a study of the properties and applications of
this proximal subdifferential, let us recall Ekeland's approximate
version of the Hopf-Rinow theorem for infinite-dimensional
Riemannian manifolds, see \cite{Ekeland2}. In our proofs we will
often use Ekeland's theorem for the cases where the complete
manifold $M$ is infinite-dimensional (so we cannot ensure the
existence of a geodesic joining any two given points of $M$).
%%%%%%%%%%%%%%%%%%
\begin{thm}[Ekeland]\label{Ekeland's approximated Hopf-Rinow theorem}
If $M$ is an infinite-dimensional Riemannian manifold which is
complete and connected then, for any given point $p$, the set
$\{q\in M : q \textrm{ {\em can be joined to $p$ by a unique
minimizing geodesic}}\}$ is dense in $M$.
\end{thm}
%%%%%%%%%%%%%%%%%%

Most of the following properties are easily translated from the
corresponding ones for $M=X$ a Hilbert space (see \cite{C})
through charts. Recall that a real-valued function $f$ defined on
a Riemannian manifold is said to be convex provided its
composition $f\circ\alpha$ with any geodesic arc $\alpha:I\to M$
is convex as a function from $I\subset \mathbb{R}$ into
$\mathbb{R}$.
%%%%%%%%%%%%%%%%%
\begin{prop}\label{basic properties}
Let $M$\ be a Riemannian manifold, $p\in M$, $f,g:M \to (-\infty
,+\infty ]$\  lower semicontinuous functions. We have
\item {i)} If $f$\ is $C^2$, then $\Pf (p)=\{ df(p)\}$.
\item{ii)} If $f$\ is convex, then $\z \in \Pf (p)$ if and only if
 $f(q)\geq f(p)+\langle\z , v\rangle$\ for every $q\in M$ and $v\in\exp_p^{-1}(q)$.
\item{iii)} If $f$\ has a local minimum at $p$, then $0\in \Pf (p)$.
\item{iv)} Every local minimum of a convex function $f$ is global.
\item{v)} If $f$\ is convex and $0\in \Pf (p)$, then $p$\ is a global minimum of $f$.
\item{vi)} $\Pf (p)+\partial _Pg(p)\subseteq \partial _P(f+g)(p)$, with equality if
$f$\ or $g$ is $C^2$.
\item{vii)} $\partial _P(cf)(p)=c\Pf (p)$, for $c>0$.
\item{viii)} If $f$\ is $K$-Lipschitz, then $\Pf (p) \subset \overline B (0,K)$.
\item{ix)} $\Pf (p)$\ is a convex subset of $TM_{p}$.
\item{x)} If $\zeta\in\partial_{P}f(p)$ and $f$ is differentiable
at $p$ then $\zeta=df(p)$.
\end{prop}
%%%%%%%%%%%%%%%%%%
\begin{proof}
All the properties but perhaps $(ii)$, $(viii)$ and $(x)$ are
easily shown to be true. Property $(viii)$ follows from the fact
that $exp_p^{-1}(.)$ is almost $1$-Lipschitz when restricted to
balls of center $0_p$ and small radius.

Let us prove $(ii)$. Let $q\in M$. Let $\gamma(t)=\exp_{p}(tv)$,
$t\in[0,1]$, which is a minimal geodesic joining $p$ and $q$. The
function $f\circ \g$ is convex and satisfies $$f(\g (t))\geq f(\g
(0))+ \langle\z , tv\rangle- \s d(\g (t)),\g (0))^2=$$ $$=f(\g
(0))+ \langle\z , t\g '(0))\rangle- \s t^2$$ for some $\s
>0$ and $t>0$ small. Hence $\langle\z ,\g '(0))\rangle \in
\partial _P(f\circ \g )(0)$, and consequently (bearing in mind that
$f\circ\gamma$ is convex on a Hilbert space) $f(\g (t))\geq f(\g
(0))+ \langle\z , t\g '(0))\rangle$, which implies $f(q)\geq
f(p)+\langle\z , v\rangle$.

To see $(x)$, note that Proposition \ref{characterization of the
proximal subdiferential through C2 functions} implies that
$\zeta\in D^{-}f(p)$, that is, $\zeta$ is a viscosity
subdifferential of $f$ at $p$ in the sense of \cite{AFL2}. Then,
since $f$ is differentiable, we have that $\zeta\in
D^{-}f(p)=D^{+}f(p)=\{df(p)\}$, so we conclude that $\zeta=df(p)$.
\end{proof}
%%%%%%%%%%%%%%%%%%
The following important result is also local, it follows from
\cite[Theorem 1.3.1]{C}.
%%%%%%%%%%%%%%%%%%
\begin{thm}[Density Theorem]\label{density theorem}
Let $M$\ be a Riemannian manifold, $p\in M$, $f:M \to (-\infty
,+\infty ]$  a lower semicontinuous function, $\e >0$. Then there
exists a point $q$ such that $d(p,q )<\e $, $f(p)-\e \leq f(q)\leq
f(p)$, and $\Pf (q)\not= \emptyset$.
\end{thm}
%%%%%%%%%%%%%%%%%%

Now we arrive to one of the main results of this section. We are
going to extend the definition and main properties of the
Moreau-Yosida regularization (see \cite{AW} for instance) to the
category of functions defined on Riemannian manifolds.
%%%%%%%%%%%%%%%%%%
\begin{thm}\label{main theorem}
Let $M$ be a complete Riemannian manifold, $f:M \to (-\infty
,+\infty ]$ a lower semicontinuous function, bounded below by a
constant $c$. Then for every $\a >0$ the function
$$f_{\a}(x)=inf_{y\in M}\{ f(y)+\a d(x,y)^2\}$$ is bounded below
by $c$, is Lipschitz on bounded sets, and satisfies $\lim_{\a \to
+\infty}f_{\a}(x)=f(x)$.

\noindent Moreover, for every $x_0\in M$ with $\partial
_Pf_{\a}(x_0)\not= \emptyset$, there is a $y_0\in M$ such
that:
\begin{enumerate}
\item[{a)}] Every minimizing sequence $\{ y_n\}$ in the
definition of $f_{\a}(x_0)$ converges to $y_0$, and consequently
the inf is a strong minimum.
\item[{b)}] There is a minimizing geodesic $\gamma$ joining $x_0$ and $y_0$.
\item[{c)}] $f_{\a}$ is differentiable at $x_0$.
\item[{d)}] $L_{x_0y_0}[df_{\a}(x_0)]\in \partial _Pf(y_0)$.
\end{enumerate}
\end{thm}
%%%%%%%%%%%%%%%%%%
\begin{proof}
It is clear that $f_{\a}(x)\geq inf_{y\in M}\{ c+\a d(x,y)^2\}=c$,
and it is easily seen that $\lim_{\a \to +\infty}f_{\a}(x)=f(x)$.
Let $A\subset M$ be a bounded set. We have that $f_{\a}(x)\leq
f(z_0)+\a d(x,z_0)^2$ for a fixed $z_0$, hence there is a positive
$m$ such that $f_{\a}(x)\leq m$ provided $x\in A$. Let us consider
$x,y\in A$ and $\e >0$, choose $z=z_{y}\in M$ such that
$f_{\a}(y)+\e \geq f(z)+\a d(y,z)^2$. We have that $d(y,z)\leq
[{1\over \a}(f_{\a}(y)+\e -c)]^{1\over 2} \leq [{1\over
\a}(m-c+\e)]^{1\over 2}:=R$. Consequently,
$$f_{\a}(x)-f_{\a}(y)\leq f_{\a}(x)-f(z)-\a d(y,z)^2+\e \leq
f(z)+\a d(x,z)^2-f(z)-\a d(y,z)^2+\e =$$ $$=\a
(d(x,z)+d(y,z))(d(x,z)-d(y,z))+\e \leq \a (2R+diamA)d(x,y)+\e.$$
By letting $\e$ go to $0$, and changing $x$ by $y$, we get that
$f_{\a}$ is Lipschitz on $A$.

\medskip

For the second part, fix $x_0\in M$, $\z \in
\partial _Pf_{\a}(x_0)$, and a sequence $(y_{n})$ such that
$f(y_{n})+\alpha d(x_{0}, y_{n}^{2})$ converges to the inf
defining $f_{\alpha}(x_{0})$. First of all let us observe that we
can always assume that for each $n$ there is a unique minimal
geodesic $\gamma_{n}:[0,1]\to M$ joining the point $y_{n}$ to
$x_{0}$. Indeed, for each couple of points $y_{n}$, $x_{0}$ we can
apply Ekeland's Theorem \ref{Ekeland's approximated Hopf-Rinow
theorem} and lower semicontinuity of $f$ to find a point $y'_{n}$
and a unique minimal geodesic $\gamma_{n}$ joining $y'_{n}$ to
$x_{0}$ in such a way that $$ d(y_{n}, y'_{n})\leq 1/n,
\,\textrm{ and }\, f(y'_{n})+\alpha d(x_{0}, y'_{n}) \leq
f(y_{n})+\alpha d(x_{0}, y_{n}). $$
Since the sequence $(y_{n})$
realizes the inf defining $f_{\alpha}(x_{0})$, so does the
sequence $(y'_{n})$. Then we can apply the argument that follows
below to the sequence $(y'_{n})$ in order to find a point $y_{0}$
with the required properties. Finally the original sequence
$(y_{n})$ must also converge to $y_{0}$ because $d(y_{n},
y'_{n})\to 0$ as $n\to\infty$. So, to save notation, we assume
$y_{n}=y'_{n}$ for each $n$.

Because $\z \in\partial _Pf_{\a}(x_0)$, there is $\s >0$ such
that, if $y$ is in a neighborhood of $x_0$, we have $$\langle\z
,exp_{x_0}^{-1}(y)\rangle\leq f_{\a}(y)-f_{\a}(x_0)+\s d(x_0,y)^2.
\eqno(*)$$ Now define $\e _n\geq 0$ by $f_{\a}(x_0)+\e
_n^2=f(y_n)+\a d(y_n,x_{0})^2$. We have $\lim_n\e _n=0$. From
$(*)$, it follows that $$\langle\z ,exp_{x_0}^{-1}(y)\rangle\leq
f(y_n)+\a d(y_n,y)^2- [f(y_n)+ \a d(y_n,x_0)^2-\e _n^2]+\s
d(x_0,y)^2=$$ $$\a d(y_n,y)^2-\a d(y_n,x_0)^2+\s d(x_0,y)^2+\e
_n^2,$$ because $f_{\a}(y)\leq f(y_n)+\a d(y_n,y)^2$.
Particularizing for $y=exp_{x_0}(\e _nv)$, with $v\in TM_{x_{0}}$,
$||v||=1$, we have $$\e _n \langle\z ,v\rangle\leq (\s +1)\e
_n^2+\a [d(y_n,y)^2-d(y_n,x_0)^2]\ \ \ (**)$$  Now let us choose
$t_n$ close enough to $1$ in order to ensure that the function
$d(.,\hat x _n)$ is differentiable at $x_0$, where $\hat x _n=\g
_n(t_n)$. Let us denote the length of $\g _{n\mid [0,t_n]}$ by
$l_n$. By using Taylor's theorem, we have that:
$$d(y_n,y)^2-d(y_n,x_0)^2\leq (d(y,\hat x _n)+l_n)^2- l(\g _n)^2
\leq $$ $$\leq (d(y,\hat x _n)+l_n)^2-(d(x_0,\hat x _n)+l_n)^2
=\f (y)-\f (x_0) =$$ $$=\f '(x_0)(\e _nv)+\f ''(exp_{x_0}^{-1}(\l
\e _nv))(\e _nv) = \f '(x_0)(\e _nv)+\f ''(x)(\e _nv),$$ where
$x=exp_{x_0}^{-1}(\l \e _nv)$, $\f (y)=(d(y,\hat x _n)+l_n)^2$,
$\f '(x_0)=2(d(x_0,\hat x _n)+l_n){\partial d(x_0,\hat x _n) \over
\partial x}$,
\par \noindent
$\f ''(x)=2[{\partial d(x,\hat x _n) \over \partial x} ]^2+
2(d(x,\hat x _n)+l_n){\partial ^2d(x,\hat x _n) \over \partial
x^2}$, and $0<\l <1$. Hence $$d(y_n,y)^2-d(y_n,x_0)^2\leq 2\e
_n(d(x_0,\hat x _n)+l_n) {\partial d(x_0,\hat x _n) \over
\partial x}(v)+2\e _n^2 [{\partial d(x,\hat x _n) \over \partial
x}(v) ]^2+$$ $$+2(d(x,\hat x _n)+l_n){\partial ^2d(x,\hat x _n)
\over \partial x^2}(\e _nv) \leq$$ $$\leq 2\e _n(d(x_0,\hat x
_n)+l_n) {\partial d(x_0,\hat x _n) \over
\partial x}(v)+ 2(d(x,\hat x _n)+l_n){\partial ^2d(x,\hat x _n)
\over \partial x^2}(\e _nv)+ 2\e _n^2\leq $$ $$\leq 2\e
_n(d(x_0,\hat x _n)+l_n) {\partial d(x_0,\hat x _n) \over
\partial x}(v)+ [{2(d(x,\hat x _n)+l_n) \over d(x,\hat x _n)} +2]\
\e _n^2, $$ since $||{\partial d(x,\hat x _n) \over \partial
x}||=1$\ and $||{\partial ^2d(x,\hat x _n) \over \partial
x^2}||={1 \over d(x,\hat x _n)}$. On the other hand, we may
firstly assume that the sequence $\{ y_n\}$ is bounded, hence so
is $\{ l_n\}$; and secondly that $\{ \hat x _n\}$ is uniformly
away from $x_0$, hence ${2(d(x,\hat x _n)+l_n) \over d(x,\hat x
_n)}+2$ is bounded by a constant $K$. Therefore
$$d(y_n,y)^2-d(y_n,x_0)^2\leq 2\e _n(d(x_0,\hat x _n)+l_n)
{\partial d(x_0,\hat x _n) \over \partial x}(v)+K \e _n^2,$$ and
from $(**)$ we get $$\e _n\langle\z -2\a (d(x_0,\hat x _n)+l_n)
{\partial d(x_0,\hat x _n) \over \partial x},v\rangle \leq (K+ \s
+1)\e _n^2.$$ This implies that $\lim_n||\z -2\a (d(x_0,\hat x
_n)+l _n){\partial d(x_0,\hat x _n) \over \partial x}||_{x_0}=0$.
>From $||{\partial d(x_0,\hat x_n) \over \partial x}||_{x_0}=1$, it
follows that $\lim_nd(x_0,\hat x _n)+\l _n ={1 \over 2 \a}||\z
||_{x_0}$\ and $\lim_n{\partial d(x_0,\hat x _n) \over \partial
x}={\z \over ||\z ||}$. Since the $\g_n$ are minimal geodesics, we
can deduce that $y_n=exp_{x_0}[-(d(x_0,\hat x _n)+\l _n){\partial
d(x_0,\hat x _n) \over \partial x}]$. Finally, the expected $y_0$\
is necessarily $y_0=exp_{x_0}(-{1 \over 2\a}\z)$. This proves (a).

\medskip

The announced geodesic in part (b) is $\g :[0,1]\to M$ defined by
$\g (t)=exp_{x_0}(-{1 \over 2\a}t\z)$, which is minimizing because
$$d(x_0,y_0)=\lim_nd(x_0,y_n)=\lim_{n} d(x_0,\hat x _n)+l _n ={1 \over
2 \a}||\z ||_{x_0}.$$

\medskip

In order to show (c), we observe that, for $y$ near $x_0$, we have
$$f_{\a}(x_0)-f_{\a}(y)\geq f(y_0)+\a d(y_0,x_0)^2-f(y_0)-\a
d(y_0,y)^2 =\a [d(y_0,x_0)^2- d(y_0,y)^2],$$ hence, using Taylor's
theorem again,
$$f_{\a}(y)\leq f_{\a}(x_0)+\a [d(y_0,y)^2-
d(y_0,x_0)^2]=f_{\a}(x_0)+\a (H(y)-H(x_0))=$$ $$=f_{\a}(x_0)+\a
DH(x_0)(exp_{x_0}^{-1}(y))+\a D^2H(x_0)(exp_{x_0}^{-1}(y))+
o(||exp_{x_0}^{-1}(y)||^2)\leq$$ $$\leq f_{\a}(x_0)+\a
DH(x_0)(exp_{x_0}^{-1}(y))+Cd(x_0,y)^2,$$ where $H(y)=(d(\hat
y_0,y)+d(\hat y_0,y_0))^2$ is $C^2$ at $x_0$, for some $\hat y_0$
lying on $\g$. This implies that $\a DH(x_0) \in
\partial ^Pf_{\a}(x_0)$ and therefore $f_{\a}$ is differentiable
at $x_0$.

\medskip

Part (d) is trivial if $x_0=y_0$. Otherwise the function $f+\a
d(x_0,.)^2=f+\a [d(x_0,\hat x_0)+d(\hat x_0,.)]^2$ attains its
minimum at $y_0$ and therefore $$0\in \partial _P(f+\a
d(x_0,.)^2)(y_0)=\partial _Pf(y_0)+ 2\a d(x_0,y_0) {\partial
d(\hat x_0,y_0) \over \partial y} $$ since $[d(x_0,\hat
x_0)+d(\hat x_0,.)]^2$ is $C^2$ at $y_0$, provided that $\hat
x_0\in \g$ and is close enough to $y_0$. Then, according to the
antisymmetry property of the partial derivatives of the distance
function (see \cite[Lemma 6.5]{AFL2}), we have
$$L_{x_0y_0}[df_{\a}(x_0)]=L_{x_0y_0}[2\a
d(x_0,y_0) {\partial d(x_0,\hat y_0) \over \partial x}]=$$
$$=-2\a d(x_0,y_0) {\partial d(\hat x_0,y_0) \over \partial y}
\in \partial _Pf(y_0).$$

\medskip

Let us observe that, as a consequence of part (c), the minimizing
geodesic joining $x_0$\ and $y_0$\ is unique.
\end{proof}
%%%%%%%%%%%%%%%%%%

Now we deduce a Borwein-Preiss variational principle for lower
semicontinuous functions defined on any complete Riemannian
manifold $M$. Let us recall that, when $M$ is
infinite-dimensional, a bounded continuous function
$f:M\to\mathbb{R}$ does not generally attain any minima. In fact,
as shown recently in \cite{AC}, the set of smooth functions with
no critical points is dense in the space of continuous functions
on $M$. Therefore, in optimization problems one has to resort to
perturbed minimization results, such as Ekeland's variational
principle (which is applicable to any complete Riemannian
manifold). Apart from Ekeland's result we have at least two other
options.

If one wants to perturb a given function with a small smooth
function which has a small derivative everywhere, in such a way
that the sum of the two functions does attain a minimum, then one
can use a Deville-Godefroy-Zizler smooth variational principle
(originally proved for Banach spaces). An extension of the DGZ
smooth variational principle is established in \cite{AFL2} for
functions defined on those Riemannian manifolds which are {\em
uniformly bumpable}.

If we wish to perturb the original function with small multiples
of squares of distance functions (so that, among other interesting
properties, we get local smoothness of the perturbing function
near the approximate minimizing point) we can use the following
Borwein-Preiss type variational principle.
%%%%%%%%%%%%%%%%%%
\begin{thm}\label{Borwein-Preiss}
Let $M$ be a complete Riemannian manifold, $f:M \to \mathbb{R}$ a lower
semicontinuous function which is bounded below, and $\e >0$. Let
$x_0\in M$ be such that $f(x_0)<inf f+\e$. Then, for every $\l
>0$ there exist $z\in B(x_0,\l )$, $y\in B(z,\l )$ with
$f(y)\leq f(x_0)$, and such that the function $\f (x)=f(x)+{\e
\over \l ^2}d(x,z)^2$ attains a strong minimum at $y$.
\end{thm}
%%%%%%%%%%%%%%%%%%
\begin{proof}
Let us consider $f_{\a}$ as in Theorem \ref{main theorem}, with
$\a ={\e \over \l ^2}$. According to the Density Theorem
\ref{density theorem}, there is a $z$ such that $d(x_0,z)\leq \l$,
$\partial _Pf_{\a}(z)\not= \emptyset$, and $f_{\a}(z)\leq
f_{\a}(x_0)\leq f(x_0)$. Hence, by Theorem \ref{main theorem},
$\f$ attains a strong minimum at a point $y_0$.

Finally, $f(y_0)+ {\e \over \l ^2}d(y_0,z)^2=f_{\a}(z)\leq f(x_0)<
inf f+\e \leq f(y_0)+\e$, hence $d(y_0,z)<\l$.
\end{proof}
%%%%%%%%%%%%%%%%%%

As an application of Theorem \ref{main theorem} we next establish
three results concerning differentiability and geometrical
properties of the distance function to a closed subset $S$ of a
Riemannian manifold $M$. These properties are probably known in
the case when $M$ is finite dimensional.

%%%%%%%%%%%%%%%%%%
\begin{thm}\label{distance to a closed set}
Let $S$ be a nonempty closed subset of $M$, and $x\in M-S$. If
$\partial _Pd_S(x)\not= \emptyset$, then $d_S$ is differentiable
at $x$. Moreover, there is an $s_0\in S$ satisfying
\item{a)} Every minimizing sequence of $d_S(x)$ converges to $s_0$.
\item{b)} $d_S(x)=d(x,s_0)$ and $d(x,s)>d_S(x)$ for every $s\in S$,
$s\neq s_{0}$.
\item{c)} There is a unique minimizing geodesic joining $x$ and $s_0$.
\end{thm}
%%%%%%%%%%%%%%%%%%
\begin{proof}
Let assume that $\xi \in \partial _Pd_S(x)$, this implies that
there is a $\s >0$ such that $$d_S(y)-d_S(x)\geq \langle\xi
,exp_x^{-1}(y)\rangle-\s d(x,y)^2$$ if $y$ is near $x$. We have
that
$$d_S^2(y)-d_S^2(x)=2d_S(x)(d_S(y)-d_S(x))+(d_S(y)-d_S(x))^2\geq
2d_S(x)(d_S(y)-d_S(x))\geq $$ $$\geq 2d_S(x)[\langle\xi ,
exp_x^{-1}(y)\rangle-\s d(x,y)^2],$$ which implies that
$2d_S(x)\xi \in
\partial _Pd_S^2(x)$. On the other hand $d_S^2(y)=inf_{z\in M}\{
I_S(z)+d(z,y)^2\}=f_{\a}(y)$ with $\a =1$, $f=I_S$, where $I_{S}$
is the indicator function of $S$, that is $I_{S}(z)=0$ if $z\in
S$, and $I_{S}(z)=\infty$ otherwise. Therefore, properties (a),
(b) and (c), which are equivalent for $d_S$ and $d_S^2$, follow
from Theorem \ref{main theorem}, as well as the fact that $d_S^2$
is differentiable at $x$. Hence $d_S$ is differentiable at $x$
because $d_S(x)>0$.
\end{proof}
%%%%%%%%%%%%%%%%%%

%%%%%%%%%%%%%%%%%%
\begin{cor}
There is a dense subset of points $x\in M-S$ such that
$d_S(x)=d(x,s_{x})$ for a unique $s_{x}\in S$ and $d_{S}$ is
differentiable at $x$.
\end{cor}
%%%%%%%%%%%%%%%%%%

%%%%%%%%%%%%%%%%%%
\begin{cor}
Let $x, x_0$ be two different points of a complete Riemannian
manifold $M$. The following statements are equivalent:
\item{i)} The function $d(.,x_0)$ is subdifferentiable at $x$.
\item{ii)} The function $d(.,x_0)$ is differentiable at $x$.
\item{\ }In both cases, there is a unique minimizing geodesic joining $x$
and $x_0$.
This condition is also equivalent provided that $M$ is finite
dimensional.
\end{cor}
%%%%%%%%%%%%%%%%%%
\begin{proof}
$(i)\implies (ii)$ and the existence of a unique minimizing
geodesic follow from Theorems \ref{main theorem} and \ref{distance
to a closed set}.

Let us assume that there is a unique minimizing geodesic joining
$x$ and $x_0$. We may prolong the geodesic up to a point $\hat x$
satisfying $d(\hat x,x_0)=d(\hat x,x)+d(x,x_0)$. In order to prove
that $d(.,x_0)$ is subdifferentiable, it is enough to see that $\f
(y)=d(\hat x,x_0)-d(y,x_0)$ is superdifferentiable at $x$, which
is a consequence of the following inequalities: $$\f (y)-\f
(x)=d(\hat x,x_0)-d(y,x_0)-d(\hat x,x)\leq d(y,\hat x)-d(\hat
x,x)\leq$$ $$\leq \langle{\partial d (x,\hat x)\over \partial
x},exp_x^{-1}(y)\rangle+\s d(x,y)^2.$$
\end{proof}
%%%%%%%%%%%%%%%%%%

\medskip

Now we turn to other topics of the theory of proximal
subdifferentials. We begin with a few local results which will be
used in some proofs (such as that of the proximal mean value
theorem and the decrease principle below). The following result
can be deduced from \cite[Theorem 1.8.3]{C}.
%%%%%%%%%%%%%%%%%%%%%
\begin{thm}[Fuzzy rule for the sum]\label{fuzzy rule for the sum}
Let $f_1, f_2:M\to(-\infty, \infty]$ be two lower semicontinuous
functions such that at least one of them is Lipschitz near $x_0$.
If $\z \in
\partial _{P} (f_1+f_2)(x_0)$ then, for every $\e >0$, there exist
$x_1$, $x_2$ and $\z _1\in \partial _Pf_1(x_1)$, $\z _2\in
\partial _Pf_2(x_2)$ such that
\item{a)} $d(x_i,x_0)<\e$\ and $|f_i(x_i)-f_i(x_0)|<\e$ for $i=1,2$.
\item{b)} $||\z -(L_{x_1x_0}(\z) _1+ L_{x_2x_0}(\z _2))||_{x_0}<\e$
\end{thm}
%%%%%%%%%%%%%%%%%%

The following theorem is also local, a consequence of the Fuzzy
Chain Rule known for functions defined on Hilbert spaces
\cite[Theorem 1.9.1, p. 59]{C}.

%%%%%%%%%%%%%%%%%%
\begin{thm}[Fuzzy chain rule]\label{fuzzy chain rule}
Let $g:N\to \mathbb{R}$ be lower semicontinuous, $F:M\to N$ be
locally Lipschitz, and assume that $g$ is Lipschitz near $F(x_0)$.
Then, for every $\z\in\partial _P(g\circ F)(x_0)$ and $\e >0$,
there are $\tilde x$, $\tilde y$ and $\eta \in \partial _Pg(\tilde
y)$ such that $d(\tilde x, x_0)<\e$, $d(\tilde y,F(x_0))<\e$,
$d(F(\tilde x),F(x_0))<\e$, and $$L_{x\tilde{x}}\z \in
\partial _P [\langle L_{\tilde{y} F(x_{0})} \left(\eta \right) ,
exp_{F(x_{0})}^{-1}\circ F(.) \rangle](\tilde x)+\e B_{TM_{\tilde
x}}.$$
\end{thm}
%%%%%%%%%%%%%%%%%%

The following result, which is local as well, relates the proximal
subdifferential $\partial_{P}f(x)$ to the viscosity
subdifferential $D^{-}f(x)$ of a function $f$ defined on a
Riemannian manifold $M$; see \cite{AFL2} for the definition of
$D^{-}f(x)$ in the manifold setting.
%%%%%%%%%%%%%%%%%%%
\begin{prop}
Let $\xi _0\in D^-f(x_0)$, $\epsilon
>0$. Then there exist $x\in B(x_0,\epsilon )$\ and $\zeta \in
\partial_{P}f(x)$\ such that $|f(x)-f(x_0)|<\epsilon$\ and $||\xi
_0-L_{xx_{0}}(\zeta )||_{x_0}$.
\end{prop}
%%%%%%%%%%%%%%%%%%
\begin{proof}
This follows from \cite[Proposition 3.4.5, p. 138]{C}.
\end{proof}
%%%%%%%%%%%%%%%%%%

Next we give a mean value theorem for the proximal subdiferential.

%%%%%%%%%%%%%%%%%%
\begin{thm}[Proximal Mean Value Theorem]\label{proximal mean value theorem}
Let $x,y\in M$, $\gamma :[0,T]\to M$ be a path joining $x$ and
$y$. Let $f$ be a Lipschitz function around $\gamma [0,T]$. Then,
for every $\varepsilon >0$, there exist $t_0$, $z\in M$ and
$\zeta \in
\partial _Pf(z)$ with
$d(z,\gamma (t_0))<\varepsilon$, and so that
$$\frac{1}{T}\left(f(y)-f(x)\right) \leq \langle \zeta ,
L_{\gamma(t_0),z}(\gamma '(t_0))\rangle +\varepsilon.$$
\end{thm}
%%%%%%%%%%%%%%%%%%
\begin{proof}
Let us consider the function $\varphi:[0,T]\to\mathbb{R}$ defined
as $$\varphi(t)=f(\gamma(t))-G(t),$$ where
$$G(t)=\frac{t}{T}f(y)+\frac{T-t}{T}f(x).$$ The function $\varphi$
is continuous, and $\varphi(0)=\varphi(T)=0$. Since the interval
$[0,T]$ is compact, there exists $t_{0}\in [0,T]$ such that
$\varphi(t_{0})\leq\varphi(t)$ for all $t\in [0,T]$. We will
consider two cases.

\noindent {\bf Case 1}. Assume that $t_{0}\in (0, T)$. Since
$\varphi$ attains a local minimum at $t_{0}$, we know that
$0\in\partial_{P}\varphi(t_{0})$. Since the function $G(t)$ is of
class $C^2$, according to the {\em easy sum rule} Proposition
\ref{basic properties}(vi), we have that
$$\frac{1}{T}(f(y)-f(x))=0+G'(t_{0})\in
\partial_{P}(f\circ\gamma)(t_{0}).$$
Now, by the fuzzy chain rule Theorem \ref{fuzzy chain rule}, there
exist $\tilde{t}$, $\tilde{z}$, $\zeta\in
\partial_{P}f(\tilde{z})$ such that
$|\tilde{t}-t_{0}|<\varepsilon$, $d(\tilde{z},
\gamma(t_{0}))<\varepsilon$,
$d(\gamma(\tilde{t}),\gamma(t_{0}))<\varepsilon$, and
\begin{eqnarray*}
& &\frac{1}{T}\left(f(y)-f(x)\right)\in
    \partial_{P}\left( \langle L_{\tilde{z}\gamma(t_{0})}(\zeta),
    \exp_{\gamma(t_{0})}^{-1}\circ\gamma(\cdot)\rangle\right)(\tilde{t})
    +[-\varepsilon, \varepsilon]=\\
& &\frac{d}{dt}\left( \langle L_{\tilde{z}\gamma(t_{0})}(\zeta),
    \exp_{\gamma(t_{0})}^{-1}\circ\gamma(\cdot)\rangle
    \right)_{|_{t=\tilde{t}}}+[-\varepsilon, \varepsilon]=\\
& &\langle L_{\tilde{z}\gamma(t_{0})}(\zeta),
\gamma'(t_{0})\rangle+[-\varepsilon, \varepsilon]= \langle \zeta,
L_{\gamma(t_{0})\tilde{z}}\left(\gamma'(t_{0})\right)\rangle
+[-\varepsilon, \varepsilon].
\end{eqnarray*}
In particular we obtain that $$\frac{1}{T}\left(f(y)-f(x)\right)
\leq \langle \zeta , L_{\gamma(t_0),z}(\gamma '(t_0))\rangle
+\varepsilon.$$

\noindent {\bf Case 2}. Now let us suppose that $t_{0}=0$ or
$t_{0}=T$. Since $\varphi(0)=\varphi(T)=0$, this means that
$\varphi(t)\geq \varphi(0)=\varphi(T)$ for all $t\in [0,T]$. We
may assume that $\varphi$ attains no local minima in $(0,T)$
(otherwise the argument of Case 1 applies and we are done). Then
there must exist $t_{0}'\in (0,T)$ such that $\varphi$ is
increasing on $(0,t_{0}')$ and $\varphi$ is decreasing on
$(t_{0}',T)$. This implies that $\zeta\geq 0$ for every
$\zeta\in\partial_{P}f(t)$ with $t\in (0,t_{0}')$, and $\eta\leq
0$ for every $\eta\in\partial_{P}f(t')$ with $t'\in (t_{0}',T)$.
Indeed, assume for instance that $t\in (0,t_{0}')$ and take
$\zeta\in\partial_{P}f(t)$. Then we have that $f(s)\geq
f(t)+\zeta(s-t)-\sigma (s-t)^{2}$ for some $\sigma\geq 0$ and $s$
in a neighborhood of $t$. By taking $s$ close enough to $t$ with
$s<t$, we get
    $$
    \zeta\geq\frac{f(t)-f(s)}{t-s}-\sigma (t-s),
    $$
and hence
    $$
    \zeta\geq\liminf_{s\to t^{-}}\left[\frac{f(t)-f(s)}{t-s}
    -\sigma (t-s)\right]\geq 0.
    $$

Now, by the Density Theorem \ref{density theorem} there exist
$t_{1}, \eta_{1}$ such that $t_{1}\in (0,t_{0}')$ and
$\eta_{1}\in\partial_{P}\varphi(t_{1})$. According to the
preceding discussion, we have $\eta_{1}\geq 0$. Since $G(t)$ is of
class $C^2$, by the {\em easy sum rule} Proposition \ref{basic
properties}(vi), we have that
$$\eta_{1}+\frac{1}{T}(f(y)-f(x))=\eta_{1}+G'(t_{1})\in
\partial_{P}(f\circ\gamma)(t_{1}).$$
Finally, by the fuzzy chain rule Theorem \ref{fuzzy chain rule},
there exist $\tilde{t}$, $\tilde{z}$, $\zeta\in
\partial_{P}f(\tilde{z})$ such that
$|\tilde{t}-t_{1}|<\varepsilon$, $d(\tilde{z},
\gamma(t_{1}))<\varepsilon$,
$d(\gamma(\tilde{t}),\gamma(t_{1}))<\varepsilon$, and
\begin{eqnarray*}
& &\eta_{1}+\frac{1}{T}\left(f(y)-f(x)\right)\in
    \partial_{P}\left( \langle L_{\tilde{z}\gamma(t_{1})}(\zeta),
    \exp_{\gamma(t_{1})}^{-1}\circ\gamma(\cdot)\rangle\right)(\tilde{t})
    +[-\varepsilon, \varepsilon]=\\
& &\frac{d}{dt}\left( \langle L_{\tilde{z}\gamma(t_{1})}(\zeta),
    \exp_{\gamma(t_{1})}^{-1}\circ\gamma(\cdot)\rangle
    \right)_{|_{t=\tilde{t}}}+[-\varepsilon, \varepsilon]=\\
& &\langle L_{\tilde{z}\gamma(t_{1})}(\zeta),
\gamma'(t_{1})\rangle+[-\varepsilon, \varepsilon]= \langle \zeta,
L_{\gamma(t_{1})\tilde{z}}\left(\gamma'(t_{1})\right)\rangle
+[-\varepsilon, \varepsilon].
\end{eqnarray*}
In particular we get $$\frac{1}{T}\left(f(y)-f(x)\right)\leq
\eta_{1}+\frac{1}{T}\left(f(y)-f(x)\right) \leq \langle \zeta ,
L_{\gamma(t_1),z}(\gamma '(t_1))\rangle +\varepsilon.$$
\end{proof}
%%%%%%%%%%%%%%%%%%

The following result is the cornerstone in the proof of the
Solvability Theorem stated below, which in turn will be the basis
of the proofs of the applications we will present later on about
circumcenters and fixed point theorems.

%%%%%%%%%%%%%%%%%%
\begin{thm}[Decrease Principle]\label{decrease principle}
Let $M$ be a complete Riemannian manifold. Let $f:M\to (-\infty
,+\infty ]$ be a lower semicontinuous function, and $x_0\in
\textrm{dom}f$. Assume that there exist $\d >0$ and $\rho >0$
satisfying that $||\zeta ||_x\geq \d$ provided that $\zeta \in
\partial _Pf(x)$, $d(x,x_0)<\rho$. Then $\inf\{ f(x):
d(x,x_0)<\rho \} \leq f(x_0)-\rho \d$.
\end{thm}
%%%%%%%%%%%%%%%%%%
\begin{proof}
We can obviously assume that $f(x_{0})=0$. Suppose that we had
$\alpha:=\inf\{f(x): d(x,x_{0})<\rho\}>-\rho\delta$. Set
$$\varepsilon=\min\{\frac{\alpha+\rho\delta}{2}, \rho\delta\}>0,
\textrm{ and }\, \varepsilon'=\frac{\varepsilon}{6(\rho+1)},$$ and
consider the function $\varphi:M\to (-\infty ,+\infty ]$ defined
by $$\varphi(x)=f(x)+(\delta-\frac{\varepsilon}{\rho})
d(x,x_{0}).$$ This function is lower semicontinuous and therefore,
by Ekeland's Variational Principle \cite{Ekeland1} it attains an
almost minimum $x_{1}$ on the closed ball
$\overline{B}(x_{0},\rho)$, that is, the function
$\varphi+\varepsilon' d(\cdot, x_{1})$ attains a minimum at
$x_{1}$ on $\overline{B}(x_{0},\rho)$. Since
    $$
    \varphi(x_{0})+\varepsilon'd(x_{0}, x_{1})\leq
    0+\varepsilon' \rho<\varepsilon\leq \alpha+\delta\rho-\varepsilon\leq
    f(y)+(\delta-\frac{\varepsilon}{\rho})\rho=\varphi(y)
    $$
for every $y$ with $d(y,x_{0})=\rho$, it is clear that $d(x_{1},
x_{0})<\rho$. Then $\varphi+\varepsilon' d(\cdot, x_{1})$ attains
a local minimum at $x_{1}$, hence
$0\in\partial_{P}(\varphi+\varepsilon' d(\cdot, x_{1}))(x_{1})$.
Now we can apply the Fuzzy Sum Rule Theorem \ref{fuzzy rule for
the sum} to obtain points $x_{2}, x_{3}$ in the open ball
$B(x_{0},\rho)$ and subdifferentials
$\zeta_{2}\in\partial_{P}\varphi(x_{2}), \,
\zeta_{3}\in\partial_{P}\varepsilon' d(\cdot, x_{1})(x_{3})$ such
that
$\|L_{x_{2}x_{1}}\zeta_{2}-L_{x_{3}x_{1}}\zeta_{3}\|<\varepsilon'$.
Since the function $\varepsilon' d(\cdot, x_{1})$ is
$\varepsilon'$-Lipschitz, it follows that
$\|\zeta_{3}\|_{x_{3}}\leq\varepsilon'$, hence
$\|\zeta_{2}\|_{x_{2}}\leq 2\varepsilon'$.

Next we consider $\zeta_{2}\in\partial_{P}\varphi(x_{2})$ and we
again apply the Fuzzy Sum Rule to the function
$\varphi=f+(\delta-\frac{\varepsilon}{\rho}) d(\cdot, x_{0})$ in
order to obtain points $x_{4}, x_{5}$ in the open ball
$B(x_{0},\rho)$ and subdifferentials
$\zeta_{4}\in\partial_{P}f(\zeta_{4}), \,
\zeta_{5}\in\partial_{P}(\delta-\frac{\varepsilon}{\rho}) d(\cdot,
x_{0})(\zeta_{5})$ such that
    $$
    \|\zeta_{2}-L_{x_{4}x_{2}}\zeta_{4}-L_{x_{5}x_{2}}\zeta_{5}\|_{x_{2}}<
    \varepsilon'.
    $$
But then, bearing in mind that $(\delta-\frac{\varepsilon}{\rho})
d(\cdot, x_{0})$ is $(\delta-\frac{\varepsilon}{\rho})$-Lipschitz,
hence $\|\zeta_{5}\|_{x_{5}}\leq
(\delta-\frac{\varepsilon}{\rho})$, it follows that
\begin{eqnarray*}
& &\|\zeta_{4}\|_{x_{4}}=\|L_{x_{4}x_{2}}\zeta_{4}\|_{x_{2}}\leq
\|\zeta_{2}\|_{x_{2}}+\|L_{x_{5}x_{2}}\zeta_{5}\|_{x_{2}}+\varepsilon'\leq
2\varepsilon'+\|\zeta_{5}\|_{x_{5}}+\varepsilon'\leq\\ & &
3\varepsilon'+(\delta-\frac{\varepsilon}{\rho})\leq
\frac{\varepsilon}{2(\rho+1)}+(\delta-\frac{\varepsilon}{\rho})\leq
\delta-\frac{\varepsilon}{2\rho}<\delta,
\end{eqnarray*}
that is, $\|\zeta_{4}\|<\delta$, a contradiction.
\end{proof}
%%%%%%%%%%%%%%%%%%

Under the same conditions we have the following Corollary (which
can also be deduced from the Variational Principle Theorem
\ref{Borwein-Preiss}).

%%%%%%%%%%%%%%%%%%
\begin{cor}
Let $\varepsilon>0$ and $x_0$ satisfy $f(x_0)<\inf f
+\varepsilon$. For every $\lambda >0$, there exist $z\in
B(x_0,\lambda )$ and $\zeta \in \partial_{P}f(z)$ such that
$f(z)<\inf f+\varepsilon$ and $\|\zeta \|<{\varepsilon/\lambda}$.
\end{cor}
%%%%%%%%%%%%%%%%%%
\begin{proof}
Otherwise, there is $\lambda>0$ so that for every $z\in
B(x_0,\lambda )$ and every $\zeta \in \partial_{P}f(z)$ we have
$\|\zeta\|\geq {\varepsilon \over \lambda}$ (we may assume that
$f(z)<\inf f+\varepsilon$\ by lower semicontinuity). Then, by the
Decrease Principle, we have $$\inf_{B(x_0,\lambda )}f(x)\leq
f(x_0)-\lambda {\varepsilon \over \lambda}<\inf f,$$ a
contradiction.
\end{proof}
%%%%%%%%%%%%%%%%%%

Now, from the Decrease Principle, we are going to obtain important
information about solvability of equations on any complete
Riemannian manifold $M$.

Let $U\subset M$, $A$ be an arbitrary set of parameters $\alpha$.
Let $F:M \times A \to [0,+\infty]$ be a function satisfying that
for every $\alpha \in A$ the function $F_{\alpha}:M\to
[0,+\infty]$ defined by $F_{\alpha}(x)= F(x,\alpha)$ is lower
semicontinuous and proper (not everywhere $\infty$). We denote the
set $\{ x\in U: F(x,\alpha )=0\}$ by $\phi (\alpha )$. Then we
have the following.

%%%%%%%%%%%%%%%%%%
\begin{thm}[Solvability Theorem]\label{solvability theorem}
Let $V\subset M$\ and $\delta >0$. Assume that $$\alpha \in A, \ \
x\in V, \ \ F(x,\alpha )>0,\ \ \zeta \in
\partial_{P}F_{\alpha}(x) \ \ \Rightarrow \ \ \|\zeta\|\geq
\delta.$$ Then for every $x\in M$\ and $\alpha \in A$, we have
$$\min\{ d(x,V^c),d(x,U^c),d(x,\phi (\alpha ))\} \leq {F(x,\alpha
)\over \delta}.$$
\end{thm}
%%%%%%%%%%%%%%%%%%
\begin{proof}
Otherwise there exist $x_0$, $\alpha _0$ and $\rho >0$ such that
$$\min\{ d(x_0,V^c),d(x_0,U^c),d(x_0,\phi (\alpha _0))\} > \rho >
{F(x_0,\alpha _0)\over \delta},$$ and in particular $B(x_0,\rho
)\subset U\cap V$ and $d(x_0,\phi (\alpha _0))>0$, which implies
that $F(x,\alpha _0)>0$ for every $x\in B(x_0,\rho )$. Hence we
have $\|\zeta\|\geq \delta$ for every $\zeta \in
\partial_P{F}_{\alpha}(x)$ with $x\in B(x_0,\rho )$. Therefore, by
the Decrease Principle, $0\leq\inf_{x\in B(x_0,\rho )}F \leq
F(x_0,\alpha _0)-\rho \delta <0$, which is a contradiction.
\end{proof}
%%%%%%%%%%%%%%%%%
Of course, the most interesting fact about the above inequality is
$d(x,\phi(\alpha))\leq F(x,\alpha)/\delta$, which implies that
$\phi(\alpha)$ is nonempty. The situation in which the above
Theorem is most often applied is when we have identified a point
$(x_{0},\alpha_{0})$ at which $F=0$, with $V$ and $\Omega$ being
neighborhoods of $x$. This is especially interesting in the cases
when the functions involved are not known to be smooth or the
derivatives do not satisfy the conditions of the implicit function
theorem. For instance, we can deduce the following result.
%%%%%%%%%%%%%%%%%
\begin{cor}
Let $x_0\in M$, $\varepsilon >0$, and $\delta >0$. Assume that
$$\alpha \in A, \ \ d(x,x_0)<2\varepsilon , \ \ F(x,\alpha )>0,\ \
\zeta \in \partial_{P}F_{\alpha}(x) \ \ \Rightarrow \ \ \|\zeta
\|\geq \delta.$$ Then we have that the equation $F(z,\alpha )=0$
has a solution for $x$ in $B(x_0,2\varepsilon )$ provided that
there is an $x\in B(x_0,\varepsilon )$ satisfying $F(x,\alpha
)<\varepsilon \delta$.
\end{cor}
%%%%%%%%%%%%%%%%%
\begin{proof}
It is enough to apply the Solvability Theorem with
$U=V=B(x_0,2\varepsilon )$. We have that $\min\{
d(x,V^c),d(x,U^c),d(x,\phi (\alpha ))\}
<\varepsilon$, and consequently $d(x,\phi (\alpha ))<\varepsilon$,
because both $d(x,U^c)$ and $d(x,V^c)$ are greater than
$\varepsilon$. Hence $\phi (\alpha )\neq \emptyset$.
\end{proof}
%%%%%%%%%%%%%%%%%
\begin{rem}
{\em Let us observe that these conditions hold if $A$ is a
topological space, $F$ is continuous, $F(x_0,\alpha _0)=0$ and
$\alpha$ is close enough to $\alpha _0$.}
\end{rem}
%%%%%%%%%%%%%%%%%

\medskip

\section{Stability of existence of circumcenters on Riemannian manifols}

\smallskip

Next we provide two applications of the Solvability Theorem. The
first one concerns some conditions on stability of the existence
of a circumcenter for three points not lying on the same geodesic.
We say that $x_{0}$ is a circumcenter for the points $a_{1},
a_{2}, a_{3}$ in a Riemannian manifold $M$ provided that $d(a_{1},
x_{0})=d(a_{2}, x_{0})=d(a_{3}, x_{0})$. It is well known that
every three points that are not aligned in $\mathbb{R}^{n}$ always
have a circumcenter. However, if $M$ is a Riemannian manifold,
this is no longer true, in general. For instance, if $M$ is the
cylinder $x^{2}+y^{2}=1$ in $\mathbb{R}^{3}$, $a_{1}=(1,0,0)$,
$a_{2}=(1,0,-1)$,
$a_{\varepsilon}=(\sqrt{1-\varepsilon^{2}},\varepsilon,1)$ and
$\varepsilon>0$ is small enough then the points $a_{1}, a_{2},
a_{\varepsilon}$ do not lie on the same geodesic and yet they have
no circumcenter.

In the sequel the symbol cut$(x_{0})$ stands for the cut locus of
the point $x_{0}$ in $M$.
%%%%%%%%%%%%%%%%%
\begin{cor}
Assume that the points $a_{1}, a_{2}, a_{3}$ have the point
$a_{0}$ as a circumcenter in a complete Riemannian manifold $M$,
and $a_{i}\notin \textrm{cut}(a_{0})$ for $i=1,2,3$. Then there is
some $\varepsilon>0$ such that for every $x_{i}\in
B(a_{i},\varepsilon)$, $i=1,2,3$, the points $x_{1}, x_{2}, x_{3}$
have a circumcenter at a point $x_{0}\in B(a_{0},\varepsilon)$.
\end{cor}
%%%%%%%%%%%%%%%%%
\begin{proof}
We may assume that $x_{1}=a_{1}$ and $x_{2}=a_{2}$. Indeed, once
the result is established in the case when one moves only one
point (for instance $a_{3}$), we can apply two times this
seemingly weaker result in order to deduce the statement of the
Corollary.

Let us pick vectors $w_{i}\in TM_{a_{0}}$ such that
$a_{i}=\exp_{a_{0}}(w_{i})$ for $i=1,2,3$, and set
$v_{i}=(1/\|w_{i}\|) w_{i}$. Since the points $a_{i}$ have the
point $a_{0}$ as a circumcenter, we have that
$\|w_{1}\|=\|w_{2}\|=\|w_{3}\|$, and the $w_{i}$ are not aligned.
Because of the strict convexity of the ball of $TM_{a_{0}}$, we
have that $v_{1}-v_{2}+\lambda(v_{1}-v_{3})\neq 0$ and
$\mu(v_{1}-v_{2})+(v_{1}-v_{3})\neq 0$ for all
$\lambda,\mu\in[-1,1]$. Then we can choose $\Delta>0$ so that
$\|v_{1}-v_{2}-\lambda(v_{1}-v_{3})\|\geq\Delta$, and
$\|\mu(v_{1}-v_{2})+(v_{1}-v_{3})\|\geq\Delta$, for every
$\lambda,\mu\in[-1,1]$. Now, by continuity of $\exp$, we can
choose numbers $\varepsilon,\delta>0$ in such a way that, if $x\in
B(a_{0},\varepsilon)$, then
\begin{enumerate}
\item[{(i)}]
$\|v_{1}(x)-v_{2}(x)-\lambda(v_{1}(x)-v_{3}(x))\|\geq \delta$ for
all $\lambda\in[-1,1]$;
\item[{(ii)}]
$\|\mu(v_{1}(x)-v_{2}(x))-(v_{1}(x)-v_{3}(x))\|\geq\delta$ for
all $\mu\in[-1,1]$,
\end{enumerate}
where $v_{i}(x)=(1/\|w_{i}(x)\|) w_{i}(x)$ and the $w_{i}(x)\in
TM_{x}$ satisfy $x_{i}=\exp_{x}(w_{i}(x))$ for $i=1,2,3$; and
$x_{1}=a_{1}, x_{2}=a_{2}$.

Let us define the function
    $$
    f(x,x_{3})=|d(x,a_{1})-d(x,a_{2})|+|d(x,a_{1})-d(x,x_{3})|,
    $$
and consider also the partial function $x\mapsto
f_{x_{3}}(x):=f(x,x_{3})$. Since the points $a_{i}$ have $a_{0}$
as a circumcenter, we have that $f(a_{0},a_{3})=0$, and, by
continuity of $f$, there is some $\eta\in (0,\varepsilon)$ such
that $f(a_{0},x_{3})<\varepsilon\delta$ for all $x_{3}\in
B(a_{3},\eta)$. Now set $V=B(a_{0},\varepsilon)$,
$A=B(a_{3},\eta)$. Then, for every $(x,x_{3})\in V\times A$ with
$f(x,x_{3})>0$ we have that $\|\zeta\|\geq\delta$ for all
$\zeta\in\partial_{P}f_{x_{3}}(x)$, because of properties $(i)$,
$(ii)$ above and the fact that
    $$
    \partial_{P}f_{x_{3}}(x)\subseteq\{\lambda(v_{1}(x)-v_{2}(x))+\mu(v_{1}(x)-v_{3}(x))
    : \max\{|\lambda|,|\mu|\}=1 \}.
    $$
Therefore, according to the Solvability Theorem \ref{solvability
theorem}, we have that
    $$
    \min\{\varepsilon, d(a_{0}, \Phi(x_{3}))\}\leq\frac{f(a_{0},x_{3})}{\delta}
    <\frac{\varepsilon\delta}{\delta}=\varepsilon,
    $$
which means that $\Phi(x_{3})\neq\emptyset$ for every $x_{3}\in
A=B(a_{0},\eta)$, that is, the points $a_{1},a_{2},x_{3}$ have a a
circumcenter at a point $x_{0}$ whose distance from $a_{0}$ is
less than $\varepsilon$, for each $x_{3}\in B(a_{0},\varepsilon)$.
\end{proof}
%%%%%%%%%%%%%%%%%
\begin{rem}
{\em The hypothesis $a_{i}\notin\textrm{cut}(a_{0})$ for $i=1,2,3$
in the above Corollary is necessary: the result may be false if
some of the $a_{i}$ belong to $\textrm{cut}(a_{0})$. For instance,
if $M$ is the cylinder defined by $x^{2}+y^{2}=1$ in
$\mathbb{R}^{3}$, $a_{0}=(0,1,0)$, $a_{1}=(0,-1,-1)$,
$a_{2}=(0,-1,1)$, and we take
$a_{3}=(a_{3}^{1},a_{3}^{2},a_{3}^{3})$ with $a_{3}^{3}>1$ and
$d(a_{0}, a_{3})=d(a_{0}, a_{2})$, then it is clear that $a_{0}$
is a circumcenter for $a_{1}, a_{2}, a_{3}$, but, for every point
$x_{3}$ lying on the geodesic going from $a_{0}$ through $a_{3}$
just past $a_{3}$, no matter how close $x_{3}$ is to $a_{3}$, the
points $a_{1}, a_{2}, x_{3}$ have no circumcenter. Notice that
$a_{1}, a_{2}\in\textrm{cut}(a_{0})=\{(x,y,z): y=-1\}$.}
\end{rem}
%%%%%%%%%%%%%%%%%
\begin{rem}
{\em We note that the above result does not follow from usual
techniques of differential calculus such as applying the implicit
function theorem.}
\end{rem}
%%%%%%%%%%%%%%%%%

%%%%%%%%%%%%%%%%%
\begin{cor}
Let $M$ be a complete Riemannian manifold such that
$\textrm{cut}(x)=\emptyset$ for all $x\in M$. Then the set
$\{(x,y,z)\in M\times M\times M: x, y, z \textrm{ have a
circumcenter}\}$ is open in $M\times M\times M$.
\end{cor}
%%%%%%%%%%%%%%%%%
Notice that there are many Riemannian manifolds $M$ satisfying the
assumptions of the above Corollary (for instance, the parabolic
cylinder defined by $y=x^{2}$ in $\mathbb{R}^{3}$).

Next we give a sufficient condition for three points in a
Riemannian manifold to have a circumcenter.
%%%%%%%%%%%%%%%%%
\begin{prop}
Consider three different points $a_{1}, a_{2}, a_{3}$ in a
Riemannian manifold $M$, and define the function
    $$
    f(x)=(d(x, a_{1})-d(x, a_{2}))^{2}+
    (d(x, a_{1})-d(x, a_{3}))^{2}.
    $$
Assume that $f$ has a global minimum at a point $a_{0}$, and
$a_{i}\notin \textrm{cut}(a_{0})$ for $i=1,2,3$. Then either the
points $a_{1}, a_{2}, a_{3}$ lie on the same geodesic or the
point $a_{0}$ is a circumcenter for $a_{1}, a_{2}, a_{3}$.
\end{prop}
%%%%%%%%%%%%%%%%%
\begin{proof}
Since $f$ attains a minimum at $a_{0}$ we have that
    $$
    \frac{1}{2}df(a_{0})=(d(a_{0},a_{1})-d(a_{0},a_{2}))(v_{1}-v_{2})+
    (d(a_{0},a_{1})-d(a_{0},a_{3}))(v_{1}-v_{3})=0, \eqno(*)
    $$
where $v_{i}=(1/\|w_{i}\|) w_{i}$, and the $w_{i}$ are such that
$\exp_{a_{0}}(w_{i})=a_{i}$ for $i=1,2,3$.

If $d(a_{0},a_{1})=d(a_{0},a_{2})=d(a_{0},a_{3})$ then $a_{0}$ is
a circumcenter for $a_{1}, a_{2}, a_{3}$ and we are done.
Otherwise, we have to show that $a_{1}, a_{2}, a_{3}$ lie on the
same geodesic. We will consider three cases.

\noindent {\bf Case 1.} $d(a_{0},a_{1})-d(a_{0},a_{2})=0$ and
$d(a_{0},a_{1})-d(a_{0},a_{3})\neq 0$. Then we have
$v_{1}=v_{3}$, which means that $a_{0}, a_{1}, a_{3}$ lie on the
same geodesic. Now we have two possibilities:

{\bf Case 1.1.} $d(a_{1},a_{2})\leq d(a_{1}, a_{3})$. Then we pick
a point $x'$ on the geodesic joining $a_{1}$ and $a_{3}$ and we
get that, if $x', a_{1}, a_{2}$ do not lie on the same geodesic,
then $f(x')=(d(x', a_{1})-d(x', a_{2}))^{2}<d(a_{1},
a_{2})^{2}\leq d(a_{1}, a_{3})^{2}=f(a_{0}),$ a contradiction.
Therefore $x', a_{1}, a_{2}$ lie on the same geodesic, hence so
do $a_{1}, a_{2}, a_{3}$.

{\bf Case 1.2.} $d(a_{1},a_{2})> d(a_{1}, a_{3})$. Now we pick a
point $x''$ on the geodesic joining $a_{1}, a_{2}$. If $x'',
a_{1}, a_{3}$ do not lie on the same geodesic then
$f(x'')=(d(a_{1}, x'')-d(a_{3}, x''))^{2}<d(a_{1},
a_{3})^{2}=f(a_{0})$, which contradicts the fact that $f$ attains
a global minimum at $a_{0}$. Hence $x'', a_{1}, a_{3}$ lie on the
same geodesic, and so do $a_{1}, a_{2}, a_{3}$.

\noindent {\bf Case 2.} $d(a_{0},a_{1})-d(a_{0},a_{2})\neq 0$ and
$d(a_{0},a_{1})-d(a_{0},a_{3})= 0$. This case is analogous to Case
2.

\noindent {\bf Case 3.} $d(a_{0},a_{2})\neq d(a_{0},a_{1})\neq
d(a_{0},a_{3})$. According to $(*)$ we have that the vectors
$v_{1}-v_{2}$ and $v_{1}-v_{3}$ linearly dependent. We may
consider two situations.

{\bf Case 3.1.} If $v_{1}-v_{3}=0$ then $v_{1}=v_{3}$, which
means that $a_{0}, a_{1}, a_{3}$ lie on the same geodesic. Then,
if we pick a point $x'$ on the geodesic joining $a_{1}, a_{2}$,
we have
\begin{eqnarray*}
& &(d(a_{1}, a_{0})-d(a_{2}, a_{0}))^{2}+
    d(a_{1}, a_{3})^{2}=f(a_{0})\leq f(x')\\
& &=(d(a_{1},a_{0})-d(a_{3},a_{0}))^{2}\leq d(a_{1}, a_{3})^{2},
\end{eqnarray*}
which implies $d(a_{1}, a_{0})=d(a_{2}, a_{0})$, a contradiction
with the standing assumptions.

{\bf Case 3.2.} If $v_{1}-v_{3}\neq 0$ we have
$v_{1}-v_{3}=\lambda (v_{1}-v_{2})$ for some $\lambda\neq 0$.
Since $\|v_{1}\|=\|v_{2}\|=\|v_{3}\|=1$ and the unit ball of
$TM_{a_{0}}$ is strictly convex, we deduce that $\lambda=1$, that
is $v_{2}=v_{3}$, which means that the points $a_{0}, a_{2},
a_{3}$ lie on the same geodesic. Let us pick a point $x'$ on the
geodesic joining $a_{1}, a_{2}$. As in the previous situation 3.1
we get that either $f(x')<f(a_{0})$ (which cannot happen) or
$d(a_{2}, a_{0})=d(a_{3}, a_{0})$, a contradiction with the
hypothesis of Case 3.
\end{proof}
%%%%%%%%%%%%%%%
\begin{rem}
{\em Note that the above proof shows that the situation considered
in Case 3 can never happen if $f$ attains a global minimum at
$a_{0}$, even in the case when all the points lie on the same
geodesic.}
\end{rem}
%%%%%%%%%%%%%%%%%

\def \R{ \bf R \rm }
\def \p{ \pi }
\def \r{ \rho }
\def \t{ \theta }
\def \l{ \lambda }
\def \e{ \varepsilon }
\def \f{ \varphi }
\def \a{ \alpha }
\def \s{ \sigma }
\def \j{ \chi }
\def \m{ \mu }
\def \N{ \bf N \rm }
\def \d{ \delta}
\def \Pf{ \partial _Pf }
\def \F{ \Phi }
\def \z{ \zeta }
\def \g{ \gamma }

\medskip

\section{Some applications to fixed point theory}

\smallskip

Now we are going to show how the Solvability Theorem allows to
deduce some interesting fixed point theorems for nonexpansive
mappings and certain perturbations of such mappings defined on
Riemannian manifolds $M$.

We first need to establish a couple of auxiliary results.
%%%%%%%%%%%%%%%%%
\begin{lem}\label{first lemma for the fixed point theorem}
Let $X, Y$ be Hilbert spaces, $F:X\to Y$ Lipschitz, and
$g:Y\to\mathbb{R}$ of class $C^2$ near $F(x_{0})$. Define
$f=g\circ F: X\to\mathbb{R}$. Then
    $$
    \partial_{P}f(x_{0})\subseteq
    \partial_{P}\left(\langle dg(F(x_{0})), F(\cdot)\rangle\right)(x_{0}).
    $$
\end{lem}
%%%%%%%%%%%%%%%%%
\begin{proof}
Take $\zeta\in \partial_{P}f(x_{0})$, that is,
    $$
    f(x)-\langle\zeta, x\rangle + \sigma\|x-x_{0}\|^{2}
    \geq f(x_{0})-\langle\zeta,x_{0}\rangle
    $$
for $x$ near $x_{0}$, with $\sigma>0$. Let $S$ be the graph of the
mapping $F$, a subset of $X\times Y$. Another way of writing the
previous inequality is the following:
    $$
    g(y)-\langle\zeta, x\rangle + \sigma\|x-x_{0}\|^{2}+I_{S}(x,y)
    \geq g(F(x_{0}))-\langle\zeta,x_{0}\rangle
    $$
for $x$ near $x_{0}$, where $I_{S}$ is the indicator function of
$S$, that is, $I_{S}(x,y)=0$ if $(x,y)\in S$, and
$I_{S}(x,y)=+\infty$ otherwise. This means that the function
    $$
    H(x,y):=g(y)-\langle\zeta, x\rangle + \sigma\|x-x_{0}\|^{2}+I_{S}(x,y):=
    h(x,y)+I_{S}(x,y)
    $$
attains a local minimum at $(x_{0}, F(x_{0}))$. Therefore
\begin{eqnarray*}
& &(0,0)\in\partial_{P}H(x_{0},
F(x_{0}))=dh(x_{0},F(x_{0}))+\partial_{P}I_{S}(x_{0},F(x_{0}))=\\
& &(-\zeta, dg(F(x_{0})))+\partial_{P}I_{S}(x_{0},F(x_{0}))=
(-\zeta, dg(F(x_{0})))+N_{S}^{P}(x_{0},F(x_{0})),
\end{eqnarray*}
where $N_{S}^{P}(x_{0},F(x_{0}))$ denotes the proximal normal
cone to $S$ at $(x_{0},F(x_{0}))$, see \cite[p.22-30]{C}. This
means that $(\zeta, -dg(F(x_{0})))\in N_{S}^{P}(x_{0},F(x_{0})),$
that is (according to \cite[Proposition 1.1.5]{C}), for some
$\sigma>0$ we have
\begin{eqnarray*}
& &\langle(\zeta, -dg(F(x_{0}))), (x, F(x))-(x_{0},
F(x_{0}))\rangle\leq
    \sigma\|(x, F(x))-(x_{0},F(x_{0}))\|^{2}=\\
& &\|(x-x_{0}\|^{2}+\|F(x)-F(x_{0}\|^{2}\leq
    \sigma(1+K)\|x-x_{0}\|^{2},
\end{eqnarray*}
where $K$ is the Lipschitz constant of $F$. Therefore,
    $$
    \langle\zeta, x-x_{0}\rangle-\sigma(1+K)\|x-x_{0}\|^{2}\leq
    \langle dg(F(x_{0}), F(x)\rangle-\langle dg(F(x_{0}), F(x_{0})\rangle,
    $$
which means that $\zeta\in\partial_{P}(\langle dg(F(x_{0})),
F(\cdot)\rangle)(x_{0})$.
\end{proof}
%%%%%%%%%%%%%%%%%

%%%%%%%%%%%%%%%%%
\begin{lem}\label{distance from x to F(x)}
Let $M$\ be a Riemannian manifold, $F:M\to M$\ Lipschitz, $x_0\in
M$ satisfying that $d(x,y)$ is $C^2$ near $(x_0,F(x_0))$,
$f(x)=d(x,F(x))$. Then $$\partial _Pf(x_0)\subset v+\partial_{P}
\langle -L_{x_0F(x_0)}v, \exp_{F(x_0)}^{-1}( F(.))
\rangle(x_0) $$ where  $v={\partial
d(x_0,F(x_0))\over
\partial x}$.
\end{lem}
%%%%%%%%%%%%%%%%%
\begin{proof}
By using the preceding Lemma, we deduce that
$$\partial_{P}f(x_0)\subset
\partial _P( \langle v,\exp_{x_0}^{-1}(.)\rangle+
\langle{\partial d(x_0,F(x_0))\over \partial
y},\exp_{F(x_0)}^{-1}(F(.))\rangle)(x_0)=$$ $$=\partial_{P}(
\langle v,\exp_{x_0}^{-1}(.)\rangle+ \langle
-L_{x_0 F(x_0)}v,\exp_{F(x_0)}^{-1}(F(.))\rangle)(x_0)=$$ $$=D(\langle
v,\exp_{x_0}^{-1}(.)\rangle)(x_0)+\partial_{P}
\langle-L_{x_0 F(x_0)}v,\exp_{F(x_0)}^{-1}(F(.))\rangle(x_0)=$$
$$v+\partial _{P}\langle -L_{x_0 F(x_0)}v, \exp_{F(x_0)}^{-1}(F(.))\rangle(x_0).$$
\end{proof}
%%%%%%%%%%%%%%%%%%
Now we consider a family of Lipschitz mappings $F_{\alpha}$,
$\alpha \in A$. Let $P_{\alpha}$ denote the set of fixed points of
$F_{\alpha}$. We are going to apply the Solvability Theorem to the
function $f(x,\alpha )=d(x,F_{\alpha}(x))$, with $U=M$. Under the
hypotheses of Lemma \ref{distance from x to F(x)} we have the
following.

%%%%%%%%%%%%%%%%%
\begin{thm}\label{main theorem on the existence of fixed points}
Let $M$ be a complete Riemannian manifold, and $F_{\alpha}:M\to
M$, $\alpha \in A$, be a family of Lipschitz mappings satisfying
the hypotheses of Lemma \ref{distance from x to F(x)} at every
point $x\in V\subset M$ with $F_{\alpha}(x)\neq x$. Assume that
there is a positive $\delta$ such that $||v_{\alpha}+\zeta||\geq \delta$
for every
$$\zeta \in
\partial_{P}\langle-v_{\alpha},(L_{F_{\alpha}(x)x}\circ exp_{F_{\alpha}(x)}^{-1}\circ
F_{\alpha})(.)\rangle(x),$$ where $x\in V$, $x\not\in P_{\alpha}$,
and $v_{\alpha}={\partial d(x,F_{\alpha}(x))\over \partial x}$.
Then we have that $$\min\{ d(x,V^c),d(x,P_{\alpha})\} \leq {d(x,
F_{\alpha}(x))\over\delta}$$ for every $x\in V$\ and $\alpha \in
A$.
\end{thm}
%%%%%%%%%%%%%%%%%
\begin{proof}
This follows immediately from Lemma \ref{distance from x to F(x)}
and Theorem \ref{solvability theorem}.
\end{proof}
%%%%%%%%%%%%%%%%%
%%%%%%%%%%%%%%%%%
\begin{rem}
{\em The condition that $d(x,y)$ is $C^2$ near $(x,F(x))$ whenever
$F(x)\neq x$ is satisfied, for instance, if $M=X$ is a Hilbert
space, or if $M$ is finite dimensional and $F(x)\not\in cut(x)$.}
\end{rem}
%%%%%%%%%%%%%%%%%

The statement of Theorem \ref{main theorem on the existence of
fixed points} might seem rather artificial at first glance but, as
the rest of the section will show, it has lots of interesting
consequences.

%%%%%%%%%%%%%%%%%
\begin{cor}
Let $M$ be a complete finite dimensional manifold or $M=X$ (a
Hilbert space). Assume that $F:M\to M$\ is Lipschitz and
$F(x)\not\in cut(x)$ for every $x\in M$. Assume also that there is
a constant $0<K<1$ such that the functions $$y\mapsto \langle
{\partial d(x,F(x))\over
\partial x}, (L^{-1}_x\circ exp_{F(x)}^{-1})(F(y)) \rangle$$ are
$K$-Lipschitz near $x$. Then $F$\ has a fixed point.
\end{cor}
%%%%%%%%%%%%%%%%%
\begin{proof}
We have that $\partial_{P}\langle-{\partial d(x,F(x))\over
\partial x}, (L^{-1}_x\circ exp_{F(x)}^{-1}\circ F)(.)\rangle(x)\subset
\overline B (0,K)$, hence $\delta =1-K$ does the work (with
$V=M$).
\end{proof}
%%%%%%%%%%%%%%%%%

The following Corollary is of course well known, but still it is
very interesting that it can be proved just by using the above
results on proximal subgradients (and without requiring any
smoothness of the distance function in $M$).
%%%%%%%%%%%%%%%%%
\begin{cor}
Let $M$ be a complete Riemannian manifold, and suppose that
$F:M\to M$ is $K$-Lipschitz, with $K<1$. Then $F$ has a unique
fixed point.
\end{cor}
%%%%%%%%%%%%%%%%%%
\begin{proof}
Uniqueness follows from the fact that $d(F(x),F(y))<d(x,y)$
whenever $x\neq y$. In order to get existence, let us observe
that, in the situation of Lemma \ref{distance from x to F(x)}, if
smoothness of the distance function fails, we may use the
following fact: $$\partial_{P}f(x_0)\subset
\partial_{L}f(x_0) \subset \bigcup _{||v||=1}[v+\partial_{L}
\langle-v,(L^{-1}_{x_0}\circ exp_{F(x_0)}^{-1}\circ
F)(.)\rangle(x_0)],$$ where $\partial_{L}f(x_0)$, the limiting
subdifferential, is defined locally in the natural way, see
\cite[p.61]{C} for the definition of $\partial_{L} f(x_{0})$ in
the Hilbert space. Next we observe that the function $$x\mapsto
\langle-v,(L^{-1}_{x_0}\circ exp_{F(x_0)}^{-1}\circ
F)(x)\rangle$$ is $(K+\varepsilon )$-Lipschitz near $x_0$, with
($K+\varepsilon
<1$), hence $$\partial_{L}\langle-v,(L^{-1}_{x_0}\circ
exp_{F(x_0)}^{-1}\circ F)(.)\rangle(x_0)\subset \overline B
(0,K+\varepsilon).$$ Therefore $\bigcup
_{||v||=1}[v+\partial_{L}\langle-v,(L^{-1}_{x_0}\circ
exp_{F(x_0)}^{-1}\circ F)(.)\rangle(x_0)]$\ does not meet the ball
$B(0,1-K-\varepsilon )$, and consequently neither does
$\partial_{P}f(x_0)$, so we can apply the Solvability Theorem as
well.
\end{proof}
%%%%%%%%%%%%%%%%%%%

\medskip

The following results, which are also consequences of Theorem
\ref{main theorem on the existence of fixed points}, allow us to
explore the behavior of small Lipschitz perturbations of certain
mappings with fixed points. Let us first observe that very small
Lipschitz perturbations of mappings having fixed points may lose
them: consider for instance $f:\mathbb{R}\to\mathbb{R}$,
$f(x)=x+\varepsilon$. The proofs of these results in their most
general (and powerful) form are rather technical. In order that
the main ideas of the proofs become apparent to the reader, we
will first establish the main theorem and its corollaries in the
case of $C^1$ smooth mappings of a Hilbert space, and then we will
proceed to study more general versions for nonsmooth perturbations
on Riemannian manifolds.

%%%%%%%%%%%%%%%%%
\begin{thm}\label{supertheorem for the Hilbert space}
Let $X$ be a Hilbert space, $G:X\to X$ a $C^1$ smooth mapping, and
$J:X\times X\to X$ satisfying that
\begin{enumerate}
\item[{(i)}] $G$ is $C$-Lipschitz on $B(x_{0},R)$;
\item[{(ii)}] $\langle h, DG(x)(h)\rangle\leq K<1$ for every $x\in
B(x_{0},R)$ and $\|h\|=1$;
\item[{(iii)}] the mapping $J_{y}:X\to X$, $J_{y}(x)=J(x,y)$ is
$L$-Lipschitz for all $y\in X$;
\item[{(iv)}] the mapping $J_x :X\to X$, $J_x(y)=J(x,y)$ is differentiable, and
 $\|\frac{\partial J}{\partial y}(x,y)-I\|\leq\varepsilon/C$ for every $x\in B(x_{0},R)$ and
$y\in G(B(x_{0},R))$;
\item[{(v)}] $L+K+\varepsilon<1$, and
\item[{(vi)}] $\|J(x_{0},G(x_{0}))-x_{0}\|< R(1-(L+K+\varepsilon))$.
\end{enumerate}
Then the mapping $F:M\to M$, defined by $F(x)=J(x,G(x))$ has a
fixed point in the ball $B(x_{0},R)$.
\end{thm}
%%%%%%%%%%%%%%%%%
\begin{proof}
This Theorem, as it is stated (that is, without assuming $J$ is
differentiable), is a consequence of Theorem \ref{supertheorem
for manifolds} below. We will make an easy proof of this
statement in the case when $J$ is differentiable and we are in a
Hilbert space setting. Let us take a
$\zeta\in\partial_{P}(\langle -v, F(\cdot)\rangle)(x)=\langle -v,
dF(x)(\cdot)\rangle$, that is,
\begin{eqnarray*}
& &\zeta=\langle -v, dF(x)(\cdot)\rangle=\langle -v,
\frac{\partial J}{\partial x}(x,G(x))(\cdot)+\frac{\partial
J}{\partial y}(x,G(x))(DG(x)(\cdot))\rangle=\\ & &\langle-v,
\frac{\partial J}{\partial x}(x,G(x))(\cdot)\rangle+\langle -v,
\frac{\partial J}{\partial y}(x,G(x))(DG(x)(\cdot))\rangle.
\end{eqnarray*}
Then, bearing in mind that $x\mapsto J_{y}(x)$ is $L$-Lipschitz
and $\partial J/\partial y$ is $\varepsilon/C$-close to the
identity, we have
\begin{eqnarray*}
& &\langle\zeta,-v\rangle= \langle-v, \frac{\partial J}{\partial
x}(x,G(x))(-v)\rangle+\langle -v, \frac{\partial J}{\partial
y}(x,G(x))(DG(x)(-v))\rangle\leq\\ &
&L+\langle-v,DG(x)(-v)\rangle+\frac{\varepsilon}{C}\|DG(x)\|\leq
L+K+\varepsilon<1.
\end{eqnarray*}
Therefore, $\|v+\zeta\|\geq\langle v, v+\zeta\rangle=\langle
v,v\rangle+\langle v,\zeta\rangle=1-\langle\zeta,-v\rangle\geq
1-(L+K+\varepsilon):=\delta>0$ and, according to Theorem \ref{main
theorem on the existence of fixed points} (here we take $A$ to be
a singleton), we get that $$\min\{ R ,d(x_0,P)\} \leq
{\|F(x_{0})-x_{0}\|\over \delta} ={ \|J(x_{0},G(x_{0}))-x_{0}\|\over
\delta}<R,$$ and consequently $P \not= \emptyset$ (that is, $F$
has a fixed point in $V=U:=B(x_{0},R)$).
\end{proof}
%%%%%%%%%%%%%%%%%
\begin{rem}\label{plain differentiability is sometimes enough}
{\em The above proof shows that Theorem \ref{supertheorem for the
Hilbert space} remains true if we only require $G$ to be
differentiable (not necessarily $C^1$) but in turn we demand that
$J$ is differentiable as well.}
\end{rem}
%%%%%%%%%%%%%%%%%
When $x_{0}$ is a fixed point of $G$ condition $(vi)$ means that
$J(x_{0}, x_{0})$ is close to $x_{0}$. The mapping $J$ can be
viewed as a general means of perturbation of the mapping $G$.
When we take a function $J$ of the form $J(x,y)=y+H(x)$ we obtain
the following Corollary, which ensures the existence of fixed
points of the mapping $G+H$ when $H$ is $L$-Lipschitz and relatively
small near $x_{0}$ (a fixed point of $G$).
%%%%%%%%%%%%%%%%%
\begin{cor}\label{main fixed point theorem for Hilbert spaces}
Let $X$ be a Hilbert space, and let $x_0$ be a fixed point of a
differentiable mapping $G:X\to X$ satisfying the following
condition: $$\langle h,DG(x)(h)\rangle\leq K<1 \ \ \ \hbox{for
every} \ \ \ x\in B(x_0,R)\ \ \hbox{and} \ \ ||h||=1.$$ Let $H$
be a differentiable $L$-Lipschitz mapping, with $L<1-K$. Then
$G+H$ has a fixed point, provided that $\|H(x_0)\|<R(1-K-L)$.
\end{cor}
%%%%%%%%%%%%%%%%%
\begin{proof}
Define $J(x,y)=y+H(x)$. Note that the above proof of Theorem
\ref{supertheorem for the Hilbert space} works for any
differentiable mappings $G$ and $J$ (not necessarily $C^1$). In
order to deduce the Corollary it is enough to observe that
$\partial H/\partial y=I$, so condition $(iv)$ of Theorem
\ref{supertheorem for the Hilbert space} is satisfied for
$\varepsilon=0$.
\end{proof}
%%%%%%%%%%%%%%%%%%

Let us observe that, when $R=+\infty$, we do not need to require
that $x_{0}$ be a fixed point of $G$, and no restriction on the
size of $H(x_0)$ is necessary either. As a consequence we have
the following.

%%%%%%%%%%%%%%%%%
\begin{cor}
Every mapping $F:X\to X$ of the form $F=T+H$, where $T$ is linear
and satisfies $\langle h,T(h)\rangle\leq K<1$\ for every
$||h||=1$, and $H$\ is $L$-Lipschitz, with $L<1-K$, has a fixed
point.
\end{cor}
%%%%%%%%%%%%%%%%%

%%%%%%%%%%%%%%%%%
\begin{rem}
{\em If $X$\ is finite dimensional, the conditions on $T$ are but
requiring that $Re\lambda <1$ for every eigenvalue $\lambda$. On
the other hand, let us observe that the function $F$\ may be
expansive, that is $||F(x)-F(y)||>||x-y||$\ for some, or even
all, $x\not= y$. Consider for instance the mapping
$T:\ell_{2}\to\ell_{2}$ defined by
    $$
    T(x_{1}, x_{2}, x_{3},x_{4}, ...)=
    5(x_{2},-x_{1},x_{4},-x_{3}, ...);
    $$
in this case $T$ is clearly expansive but we can take $K=0$. This
result should be compared with \cite[Corollary 1.6, p.24]{Dug}.}
\end{rem}
%%%%%%%%%%%%%%%%%

As a consequence of Corollary \ref{main fixed point theorem for
Hilbert spaces} we can also deduce the following local version of
the result.
%%%%%%%%%%%%%%%%%
\begin{cor}\label{local fixed point theorem}
Let $x_0$ be a fixed point of a differentiable function $G:X\to X$
satisfying the following condition: $$\langle h,DG(x_0)(h)\rangle
\leq K<1 \ \ \ \hbox{for every} \ \ \ \|h\|=1.$$ Let $H$\ be a differentiable
$L$-Lipschitz function. Then there exists a positive constant
$\alpha _0$ such that the function $G+\alpha H$ has a fixed
point, for every $\alpha \in (0,\alpha _0)$.
\end{cor}
%%%%%%%%%%%%%%%%%

\medskip

Another way to perturb a mapping $G$ with a fixed point $x_{0}$
is to compose it with a function $H$ which is close to the
identity. When we take $J$ of the form $J(x,y)=H(y)$ in Theorem
\ref{supertheorem for the Hilbert space} we obtain the following.

%%%%%%%%%%%%%%%%
\begin{cor}
Let $X$ be a Hilbert space, and $x_{0}$ be a fixed point of a
differentiable mapping $G:X\to X$ such that
    $$
    \langle h, DG(x)(h)\rangle\leq K<1
    $$
for every $x\in B(x_{0}, R)$ and $h\in X$ with $\|h\| =1$. Let
$H:X\to X$ be a differentiable mapping such that
$\|DH(G(x))-I\|<\varepsilon$ for every $x\in B(x_{0}, R)$. Then
$F=H\circ G$ has a fixed point in $B(x_{0}, R)$, provided that
$K+\varepsilon<1$ and $\|H(x_{0})-x_0\|<R(1-K-\varepsilon)$.
\end{cor}
%%%%%%%%%%%%%%%%
\begin{proof}
For $J(x,y)=H(y)$ we have that $x\mapsto J_{y}(x)$ is constant,
hence $0$-Lipschitz for every $y$, and we can apply Theorem
\ref{supertheorem for the Hilbert space} with $L=0$ (and bearing
in mind Remark \ref{plain differentiability is sometimes enough}).
\end{proof}
%%%%%%%%%%%%%%%%

\medskip

Finally we will consider an extension of Theorem \ref{supertheorem
for the Hilbert space} and the above corollaries to the setting of
nonsmooth mappings on Riemannian manifolds. We will make use of
the following fact about {\em partial proximal subdifferentials}.
%%%%%%%%%%%%%%%%
\begin{lem}\label{partial proximal subgradients}
Let $M$ be a Riemannian manifold, and $f:M\times M\to\mathbb{R}$.
For each $x\in M$ let us define the partial function
$f_{x}:M\to\mathbb{R}$ by $f_{x}(y)=f(x,y)$, and define also
$f_{y}:M\to\mathbb{R}$ by $f_{y}(x)=f(x,y)$. Assume that
$\zeta=(\zeta_{1},\zeta_{2})\in\partial_{P}f(x_{0},y_{0})$. Then
$\zeta_{1}\in\partial_{P}f_{y_{0}}(x_{0})$, and
$\zeta_{2}\in\partial_{P}f_{x_{0}}(y_{0})$.
\end{lem}
%%%%%%%%%%%%%%%%
\begin{proof}
Since $\zeta\in\partial_{P}f(x_{0},y_{0})$ there exists a $C^2$
function $\varphi:M\times M\to\mathbb{R}$ such that $f-\varphi$
attains a local minimum at $(x_{0},y_{0})$, and
$$(\frac{\partial\varphi}{\partial x}(x_{0},y_{0}),
\frac{\partial\varphi}{\partial
y}(x_{0},y_{0}))=d\varphi(x_{0},y_{0})
=\zeta=(\zeta_{1},\zeta_{2}).$$ Then it is obvious that $x\mapsto
f_{y_{0}}(x)-\varphi(x,y_{0})$ attains a local minimum at $x_{0}$
as well, so
    $$
    \zeta_{1}=\frac{\partial\varphi}{\partial x}(x_{0},y_{0})
    \in\partial_{P}f_{y_{0}}(x_{0}).
    $$
In the same way we see that
$\zeta_{2}=\frac{\partial\varphi}{\partial y}(x_{0},y_{0})\in
\partial_{P}f_{x_{0}}(y_{0})$.
\end{proof}
%%%%%%%%%%%%%%%%
Now we can prove our main result about perturbation of mappings
with fixed points. As said before, the mapping $J$ should be
regarded as a general form of perturbation of $G$. We use the following notation:
$$\textrm{sing}(x):=\{ y\in M : d(\cdot, x)^2\, \textrm{ is not differentiable at } y\}.$$
When $M$ is finite-dimensional it is well known that $\textrm{sing}(x)\subseteq
\textrm{cut}(x)$, and both $\textrm{sing}(x)$ and $\textrm{cut}(x)$ are sets of measure zero.

%%%%%%%%%%%%%%%%
\begin{thm}\label{supertheorem for manifolds}
Let $M$ be a complete Riemannian manifold, $G:M\to M$, $J:M\times
M\to M$ and $F:M\to M$ defined by $F(x)=J(x,G(x))$ be mappings
such that:
\begin{enumerate}
\item[{(i)}] $F(x)\notin\textrm{sing}(x)\cup\textrm{sing}(G(x))$ for every $x\in B(x_{0},R)$, and moreover there are (unique) minimizing geodesics joining  $F(x)$ to $x$ and $F(x)$ to $G(x)$;
\item[{(ii)}] $G$ is $C^1$ smooth and
$\langle L_{xF(x)}h, L_{G(x)F(x)}dG(x)(h)\rangle_{F(x)}\leq K<1$
for all $x\in B(x_{0},R)$ and $h\in TM_{x}$ with $\|h\|_{x}=1$;
\item[{(iii)}] $G$ is $C$-Lipschitz on $B(x_{0},R)$;
\item[{(iv)}] $J$ is locally Lipschitz;
\item[{(v)}] the mapping $x\mapsto J_{y}(x):=J(x,y)$ is $L$-Lipschitz for
every $y\in M$;
\item[{(vi)}] the mapping $y\mapsto J_{x}(y):=J(x,y)$ is differentiable for every
$x\in B(x_{0},R)$, there is a unique minimizing geodesic joining $J(x,y)$ to $y$,  $J(x,y)\notin\textrm{sing}(y)$ and
    $$
    \|\frac{\partial J}{\partial
    y}(x,y)-L_{yJ(x,y)}\|\leq\frac{\varepsilon}{C}
    $$
for every $x\in B(x_{0},R)$ and $y\in G(B(x_{0},R))$;
\item[{(vii)}] $L+K+\varepsilon<1$; and
\item[{(viii)}] $d(x_{0}, J(x_{0},G(x_{0})))<R(1-(L+K+\varepsilon))$.
\end{enumerate}
Then the mapping $F:M\to M$ has a
fixed point in the ball $B(x_{0},R)$.

Moreover, when $M$ is finite-dimensional, conditions $(i)$ and $(vi)$ on the singular sets
can be replaced by:
\begin{enumerate}
\item[{(i)'}] $F(x)\notin\textrm{cut}(x)\cup\textrm{cut}(G(x))$ for every $x\in B(x_{0},R)$, and
\item[{(vi)'}] the mapping $y\mapsto J_{x}(y):=J(x,y)$ is differentiable for every
$x\in B(x_{0},R)$, and $J(x,y)\notin\textrm{cut}(y)$ and
    $$
    \|\frac{\partial J}{\partial
    y}(x,y)-L_{yJ(x,y)}\|\leq\frac{\varepsilon}{C}
    $$
for every $x\in B(x_{0},R)$ and $y\in G(B(x_{0},R))$.
\end{enumerate}
\end{thm}
%%%%%%%%%%%%%%%%

At first glance this statement may seem to be overburdened with
assumptions, but it turns out that all of them are either useful
or necessary, as we will see from its corollaries and in the next
remarks. Before giving the proof of the Theorem, let us make some
comments on these assumptions.

\noindent {\bf 1.} The condition in $(i)'$ that
$F(x)\notin\textrm{cut}(x)$ for all $x\in B(x_{0},R)$ is
necessary. Indeed, in the simplest case when there is no
perturbation at all, that is, $J(x,y)=y$, if we take $G$ to be a
continuous mapping form the sphere $S^2$ into itself and
$G(x)\in\textrm{sing}(x)=\textrm{cut}(x)=\{-x\}$ for all $x\in S^2$, then $G$ is
the antipodal map and  has no fixed point. Therefore, in order
that $G:S^2\to S^2$ has a fixed point, there must exist some
$x_{0}$ with $G(x_{0})\notin\textrm{cut}(x_{0})$ and, therefore,
by continuity, $G(x)\notin\textrm{cut}(x)$ for every $x$ in a
neighborhood of $x_{0}$.

On the other hand, not only is this a necessary condition, but
also very natural in these kinds of problems. For instance, one
can deduce from the hairy ball theorem that if $G:S^{2}\to S^{2}$
is a continuous mapping such that $G(x)\notin\textrm{cut}(x)$ for
every $x\in S^{2}$ then $G$ has a fixed point. Indeed, for every
$x\in S^{2}$, the condition $G(x)\notin\textrm{cut}(x)$ implies
the existence of a unique $v_{x}\in TS^{2}_{x}$ with $\|v_{x}\|<\pi$ such that
$\exp_{x}(v_{x})=G(x)$. The mapping $S^{2}\ni x\mapsto v_{x}\in
TS^{2}$ defines a continuous field of tangent vectors to $S^{2}$.
If $G$ did not have any fixed point then we would have $v_{x}\neq
0$ for all $x$, which contradicts the hairy ball theorem.

\noindent {\bf 2.} The other condition in $(i)'$ that
$F(x)\notin\textrm{cut}(x)$ is also natural in this setting and
verily easily satisfied if we mean $F$ to be a relatively small
perturbation of $G$. For instance, if $M$ has a positive
injectivity radius $\rho=i(M)>0$ and $F$ is $\rho$-close to $G$,
that is, $d(F(x),G(x))<\rho$ for $x\in B(x_{0},R)$, then
$F(x)\notin\textrm{cut}(G(x))$ for $x\in B(x_{0},R)$.

\noindent {\bf 3.} Since the main aim of the present theorem is to
establish corolaries in which we have a mapping $G:M\to M$ with a
fixed point $x_{0}$ and we perturb G by composing with or summing
a mapping $H$ with certain properties, thus obtaining a mapping
$F$ which is relatively close to $G$, and then we want to be able
to guarantee that this perturbation of $G$ still has a fixed
point, it turns out that condition $(i)'$ of the Theorem is not
really restrictive. Indeed, since $G(x_{0})=x_{0}$, $J(x_{0},
x_{0})$ is relatively close to $x_{0}$ (see property $(viii)$),
the mappings $G$ and $F=J\circ G$ are continuous and there always
exists a convex neighborhood of $x_{0}$ in $M$, it is clear that
there must be some $R>0$ such that $F(x)\notin\textrm{cut}(x)$
for every $x\in B(x_{0}, R)$.

\noindent {\bf 4.} The second part of condition $(v)$ means that,
in its second variable, $J$ is relatively close to the identity, a
natural condition to put if we mean the function $F(x)=J(x,G(x))$
to be a relatively small perturbation of $G$.

\noindent {\bf 5.} The requirement that $G$ is Lipschitz on
$B(x_{0},R)$ is not a strong one. On the one hand, since $G$ is
$C^1$ it is locally Lipschitz, so condition $(iii)$ is met
provided $R$ is small enough. On the other hand, when $M$ is
finite dimensional, by local compactness of $M$ and continuity of
$dG$, condition $(iii)$ is always true for any $R$.

\noindent {\bf 6.} Condition $(ii)$ is met in many interesting
situations: for example, when the behavior of $G$ in a
neighborhood of $x_{0}$ is similar to a multiple of a rotation.
Consider for instance $M$ the surface $z=x^{2}+y^{2}$ in
$\mathbb{R}^{3}$, and $G(x,y,z)=(5y,-5x,25z)$. Then $dG(0)$ is the
linear mapping $T(x,y)=5(y,-x)$, and it is clear that for every
$K\in (0,1)$ there is some $R>0$ such that $(ii)$ is satisfied.
Of course the origin is a fixed point of $G$. Theorem
\ref{supertheorem for manifolds} tells us that any relatively
small perturbation of $G$ still has a fixed point (relatively
close to the origin).

\noindent {\bf 7.} Notice also that Theorem \ref{supertheorem for
manifolds} gives, in the case of $M=X$ a Hilbert space, the
statement of Theorem \ref{supertheorem for the Hilbert space},
which is stronger than the version already proved, because here
the mapping $J$ is not necessarily differentiable, only
${\frac{\partial J}{\partial y}}$ needs to exist. This is one of
the reasons why the proof of Theorem \ref{supertheorem for
manifolds} is much more complicated than the one already given in
the Hilbert space. This seemingly little difference is worth the
effort of the proof, because, for instance, in Theorem \ref{fixed
point theorem for Lie groups}, we only have to ask that the
perturbing function $H$ is Lipschitz, not necessarily everywhere
differentiable with a bounded derivative.
%%%%%%%%%%%%%%%%

\begin{center}
{\bf Proof of Theorem \ref{supertheorem for manifolds}}
\end{center}

Let us denote $\tilde{G}(x)=(x,G(x))$. Let us fix a point $x\in
B(x_{0},R)$ with a subgradient
    $$
    \zeta\in\partial_{P}\big(\langle-v, L_{F(x)x}\circ\exp_{F(x)}^{-1}\circ
    F(\cdot)\rangle\big)(x)=
    \partial_{P}\Big(\langle -v, L_{F(x)x}\circ\exp_{F(x)}^{-1}\circ
    J(\cdot)\rangle \circ\tilde{G}(\cdot)\Big)(x),
    $$
where $v=\partial d(x, F(x))/\partial x$. We want to see that
$\|v+\zeta\|\geq\lambda>0$ for some $\lambda>0$ independent of
$x, \zeta$, and such that $d(J(x_{0},G(x_{0})), x_{0})<R\lambda$.
Then, by using Theorem \ref{main theorem on the existence of
fixed points}, we will get that
    $$
    \min\{R, d(x_{0}, P)\}\leq\frac{d(F(x_{0}), x_{0})}{\lambda}=
    \frac{d(J(x_{0},G(x_{0})), x_{0})}{\lambda}<R,
    $$
hence $P\neq\emptyset$, that is, $F$ has a fixed point in
$V=U=B(x_{0},R)$. So let us prove that there exists such a number
$\lambda$.

Since $J$ is locally Lipschitz we can find positive numbers $C'$,
$\delta_{0}$ such that $J$ is $C'$-Lipschitz on the ball
$B(\tilde{G}(x), \delta_{0})$.

By continuity of $\tilde{G}$ and $d\tilde{G}$ and the properties
of $\exp$, for any given $\varepsilon'>0$ we can find
$\delta_{1}>0$ such that
    $$
    \|d(\exp_{\tilde{G}(x)}^{-1})(\tilde{G}(\tilde{x}))
    -L_{\tilde{G}(\tilde{x})\tilde{G}(x)}\|
    <\frac{\varepsilon'}{1+\|d\tilde{G}(x)\|},\,
     \textrm{ and }\, \|d\tilde{G}(\tilde{x})\|\leq
     1+\|d\tilde{G}(x)\|
    $$
whenever $\tilde{x}\in B(x,\delta_{1})$, and therefore, by the
chain rule,
    $$
    \|d(\exp_{\tilde{G}(x)}^{-1}\circ\tilde{G})(\tilde{x})
    -L_{\tilde{G} (\tilde{x} )\tilde{G}(x)}d\tilde{G}(\tilde{x})\|
    <\frac{\varepsilon'}{1+\|d\tilde{G}(x)\|}
    \|d\tilde{G}(\tilde{x})\|\leq\varepsilon'    \eqno(1)
    $$
for every $\tilde{x}\in B(x_{0},\delta_{1})$. On the other hand,
since $\exp_{F(x)}^{-1}$ is an almost isometry near $F(x)$ and the
mapping $J$ is continuous, we can find $\delta_{2}>0$ such that if
$\tilde{y},\tilde{z}\in B(\tilde{G}(x),\delta_{2})$ then
    $$
    \|\exp_{F(x)}^{-1}(J(\tilde{y}))-
    \exp_{F(x)}^{-1}(J(\tilde{z}))\|\leq (1+\varepsilon')
    d(J(\tilde{y}),J(\tilde{z})).
    $$
In particular, bearing in mind the fact that the mapping $J_{y}$
is $L$-Lipschitz for all $y$, we deduce that
    $$
    \|\exp_{F(x)}^{-1}(J(z,y))-
    \exp_{F(x)}^{-1}(J(z',y))\|\leq (1+\varepsilon')\,L\,
    d(z,z') \eqno(2)
    $$
for all $(z,y), (z',y)\in B(\tilde{G}(x),\delta_{2})$.

In a similar manner, because $d\exp_{F(x)}^{-1}(\tilde{z})$ is
arbitrarily close to $L_{\tilde{z}F(x)}$ provided $\tilde{z}$ is
close enough to $F(x)$, and $J(\tilde{y})$ is close to
$F(x)=J(\tilde{G}(x))$ when $\tilde{y}$ is close to
$\tilde{G}(x)$, we can find a number $\delta_{3}>0$ such that
    $$
    \|d\exp_{F(x)}^{-1}(J(\tilde{y}))-L_{J(\tilde{y})F(x)}\|
    \leq\varepsilon' \eqno(3)
    $$
provided that $\tilde{y}\in B(\tilde{G}(x),\delta_{3})$.

Because of the continuity properties of the parallel transport
and the geodesics (a consequence of their being solutions of
differential equations which exhibit continuous dependence with
respect the initial data), we may find numbers $\delta_{4},
\delta_{5}, \delta_{6}>0$ such that: $$
\|L_{J(\tilde{y})F(x)}L_{\tilde{y_2}J(\tilde{y})}-L_{\tilde{y_2}F(x)}\|
\leq\varepsilon'  \textrm{ provided that }
d(\tilde{y},\tilde{G}(x))<\delta_{4}; \eqno(4)$$ $$
\|dG(\tilde{x} )L_{x\tilde{x}}-L_{G(x)G(\tilde{x} )}dG(x)\|
\leq\varepsilon'  \textrm{ provided that } d(x,\tilde{x}
)<\delta_{5}, \textrm{ and } \eqno(5) $$ $$
\|L_{\tilde{y_2}F(x)}L_{G(x)\tilde{y_2}}-L_{G(x)F(x)}\|
\leq\varepsilon' \textrm{ provided that }
d(\tilde{y_2},G(x))<\delta_{6}. \eqno(6) $$

Let us take any $\delta<\min\{\delta_{0},
\delta_{1},\delta_{2},\delta_{3},\delta_{4},\delta_{5},\delta_{6}\}$.
By the fuzzy chain rule Theorem \ref{fuzzy chain rule}, we have
that there are points $\tilde{y}=(\tilde{y}_{1},\tilde{y}_{2})\in
M\times M$, $\tilde{x}\in M$, and a subgradient
    $$
    \eta\in\partial_{P}\big(\langle-v, L_{F(x)x}\circ\exp_{F(x)}^{-1}
    \circ J(\cdot)\rangle\big)(\tilde{y}) \eqno(7)
    $$
such that $d(\tilde{y},\tilde{G}(x))<\delta$,
$d(\tilde{x},x)<\delta$, $d(\tilde{G}(\tilde{x}),
\tilde{G}(x))<\delta$, and
    $$
    L_{x\tilde{x}} \zeta\in\partial\Big(\langle L_{\tilde{y}\tilde{G}(x)}(\eta),
    \exp_{\tilde{G}(x)}^{-1}\circ\tilde{G}\rangle\Big)(\tilde{x})
    +\delta B_{TM_{\tilde{x}}}.
    $$
Since the mapping $\tilde{G}$ is differentiable this means,
according to property $(x)$ of Proposition \ref{basic properties},
that
    $$
    L_{x\tilde{x}} \zeta \in \langle L_{\tilde{y}\tilde{G}(x)}(\eta),
    d\big(\exp_{\tilde{G}(x)}^{-1}\circ\tilde{G}\big)(\tilde{x})
    (\cdot)\rangle+\delta B_{TM_{\tilde{x}}}.
    $$
But, by equation $(1)$ above, and taking into account that
$\|\eta\|\leq C'$ (because $J$ is $C'$-Lipschitz on
$B(\tilde{G}(x),\delta)\ni\tilde{y}$), we have that
\begin{eqnarray*}
& &\langle L_{\tilde{y}\tilde{G}(x)}(\eta),
    d\big(\exp_{\tilde{G}(x)}^{-1}\circ\tilde{G}(\cdot)\big)(\tilde{x})
    (\cdot)\rangle+\delta B_{TM_{\tilde{x}}}\\
& &\subseteq \langle L_{\tilde{y}\tilde{G}(x)}(\eta),
    L_{\tilde{G}(\tilde{x})\tilde{G}(x)}\circ
    d\tilde{G}(\tilde{x})(\cdot)\rangle+(\delta+\varepsilon' \| \eta \|)
    B_{TM_{\tilde{x}}}\\
& &\subseteq \langle L_{\tilde{y}\tilde{G}(x)}(\eta),
    L_{\tilde{G}(\tilde{x})\tilde{G}(x)}\circ
    d\tilde{G}(\tilde{x})(\cdot)\rangle+(\delta+\varepsilon' C')
    B_{TM_{\tilde{x}}},
\end{eqnarray*}
so we get that, denoting $\eta=(\eta_{1},\eta_{2})$,
\begin{eqnarray*}
& &L_{x\tilde{x}} \zeta\in\langle L_{\tilde{y}\tilde{G}(x)}(\eta),
    L_{\tilde{G}(\tilde{x})\tilde{G}(x)}\circ
    d\tilde{G}(\tilde{x})(\cdot)\rangle+(\delta+\varepsilon' C')
    B_{TM_{\tilde{x}}}=\\
& &\langle (\eta_{1},\eta_{2}),
L_{\tilde{G}(x)\tilde{y}}L_{\tilde{G}(\tilde{x})\tilde{G}(x)}\circ
d\tilde{G}(\tilde{x})(\cdot)\rangle+ (\delta+\varepsilon'
C')B_{TM_{\tilde{x}}}=\\ & &\langle\eta_{1},
L_{x\tilde{y}_{1}}L_{\tilde{x}x}(\cdot) \rangle + \langle\eta_{2},
L_{G(x)\tilde{y}_{2}}L_{G(\tilde{x})G(x)}\circ
dG(\tilde{x})(\cdot)\rangle+(\delta+\varepsilon' C'
)B_{TM_{\tilde{x}}}. \hspace{0.5cm} (8)
\end{eqnarray*}
Now, from inequality $(2)$ above and taking into account that
$L_{F(x)x}$ is an isometry, we get that the functions
    $$z\mapsto\langle-v,
    L_{F(x)x}\circ\exp_{F(x)}^{-1}
    \circ J(z,y)\rangle$$
are $L(1+\varepsilon')$-Lipschitz on $B(x,\delta_2)$
for every $y\in B(G(x),\delta_2)$. Then, since
$$(\eta_{1},\eta_{2})=\eta\in\partial_{P}\big(\langle-v,
    L_{F(x)x}\circ\exp_{F(x)}^{-1}
    \circ J(\cdot)\rangle\big)(\tilde{y}), \eqno(9)
    $$
and by using Lemma \ref{partial proximal subgradients}, we deduce
that $$\|\eta_{1}\|\leq (1+\varepsilon')L. \eqno (10)$$

On the other hand, since the mapping $y\mapsto J(x,y)$ is
differentiable, by looking at equation $(9)$ above, and again
using Proposition \ref{basic properties}(x) and Lemma \ref{partial
proximal subgradients}, we have that
\begin{eqnarray*}
& &\eta_{2}=\frac{\partial\big(\langle-v,
    L_{F(x)x}\circ\exp_{F(x)}^{-1}
    \circ J(y_{1},y_{2})\rangle\big)}
    {\partial y_{2}}(\tilde{y})=\\
& &\frac{\partial\big(\langle-L_{xF(x)}v,
    \circ\exp_{F(x)}^{-1}
    \circ J(y_{1},y_{2})\rangle\big)}
    {\partial y_{2}}(\tilde{y})=\\
& &\langle -L_{xF(x)}(v),
d\exp_{F(x)}^{-1}(J(\tilde{y}))(\frac{\partial J}{\partial
y_{2}}(\tilde{y}))(\cdot)\rangle. \hspace{4.5cm} (11)
\end{eqnarray*}
Besides, bearing in mind equation $(3)$, $(4)$ and the assumption
$(vi)$ of the statement, we have
\begin{eqnarray*}
&\|d\exp_{F(x)}^{-1}(J(\tilde{y}))\circ\frac{\partial J}{\partial
y_{2}}(\tilde{y}) - L_{\tilde{y}_{2}F(x)}\|\leq \\ &
\|d\exp_{F(x)}^{-1}(J(\tilde{y}))\circ\frac{\partial J}{\partial
y_{2}}(\tilde{y}) - L_{J(\tilde{y})F(x)}\circ
L_{\tilde{y}_{2}J(\tilde{y})}\| +\varepsilon'=\\ &
\|d\exp_{F(x)}^{-1}(J(\tilde{y}))\circ\big(\frac{\partial
J}{\partial y_{2}}(\tilde{y})-L_{\tilde{y}_{2}J(\tilde{y})}\big)+
\big(d\exp_{F(x)}^{-1}(J(\tilde{y}))-L_{J(\tilde{y})F(x)}\big)\circ
L_{\tilde{y}_{2}J(\tilde{y})}\| +\varepsilon' \leq \\ &
\|d\exp_{F(x)}^{-1}(J(\tilde{y}))\|\, \|\frac{\partial
J}{\partial y_{2}}(\tilde{y})-L_{\tilde{y}_{2}J(\tilde{y})}\|+
\|d\exp_{F(x)}^{-1}(J(\tilde{y}))-L_{J(\tilde{y})F(x)}\|\,
\|L_{\tilde{y}_{2}J(\tilde{y})}\| +\varepsilon' \leq\\ &
(1+\varepsilon')\cdot\frac{\varepsilon}{C}+\varepsilon'\cdot
1+\varepsilon' =
(1+\varepsilon')\frac{\varepsilon}{C}+2\varepsilon',
\end{eqnarray*}
which, combined with $(11)$, yields
\begin{eqnarray*}
& &\langle\eta_{2}, h\rangle=\langle-L_{xF(x)}(v),
    d\exp_{F(x)}^{-1}(J(\tilde{y}))(\frac{\partial J}{\partial y_{2}}
    (\tilde{y}))(h)\rangle\leq\\
& &\langle -L_{xF(x)}(v), L_{\tilde{y_2}F(x)} h\rangle + \big(
(1+\varepsilon')\frac{\varepsilon}{C}+2\varepsilon'\big)\|h\|
\end{eqnarray*}
for all $h\in TM_{\tilde{y}_{2}}$. By taking
$h=L_{G(x)\tilde{y}_{2}} L_{G(\tilde{x})G(x)}
dG(\tilde{x})(-L_{x\tilde{x}} v)$ into this expression we get
\begin{eqnarray*}
& &\langle\eta_{2}, L_{G(x)\tilde{y}_{2}} L_{G(\tilde{x})G(x)}
dG(\tilde{x})(-L_{x\tilde{x}} v)\rangle\leq\\ & & \langle
-L_{xF(x)} v,
L_{\tilde{y_2}F(x)}L_{G(x)\tilde{y}_{2}}L_{G(\tilde{x})G(x)}
dG(\tilde{x})(-L_{x\tilde{x}}
v)\rangle+\big((1+\varepsilon')\frac{\varepsilon}{C}+
2\varepsilon'\big)\|dG(\tilde{x})\|\leq \\ & & \langle -L_{xF(x)}
v, L_{\tilde{y_2}F(x)}L_{G(x)\tilde{y}_{2}}L_{G(\tilde{x})G(x)}
dG(\tilde{x})(-L_{x\tilde{x}}
v)\rangle+\big((1+\varepsilon')\frac{\varepsilon}{C}+
2\varepsilon'\big)C =\\ & & \langle -L_{xF(x)} v,
L_{\tilde{y_2}F(x)}L_{G(x)\tilde{y}_{2}}L_{G(\tilde{x})G(x)}
dG(\tilde{x})(-L_{x\tilde{x}}
v)\rangle+(1+\varepsilon')\varepsilon+ 2\varepsilon'C.
\hspace{1.5cm}(12)
\end{eqnarray*}
Now, by combining  equations $(8)$, $(10)$, $(12)$, $(5)$, $(6)$
and assumption $(ii)$, we obtain
\begin{eqnarray*}
& &\langle\zeta,-v\rangle=\langle L_{x\tilde{x}} \zeta,
-L_{x\tilde{x}} v\rangle\leq\\ & & (\delta+\varepsilon' C')+
\langle\eta_{1},
L_{x\tilde{y}_{1}}L_{\tilde{x}x}(-L_{x\tilde{x}}v) \rangle +
\langle\eta_{2}, L_{G(x)\tilde{y}_{2}}L_{G(\tilde{x})G(x)}\circ
dG(\tilde{x})(-L_{x\tilde{x}}v)\rangle \leq\\ & &
(\delta+\varepsilon' C' ) + (1+\varepsilon') L+\langle-L_{xF(x)}v,
L_{\tilde{y_2}F(x)}L_{G(x)\tilde{y}_{2}}L_{G(\tilde{x})G(x)}\circ
dG(\tilde{x})(-L_{x\tilde{x}}v)\rangle\\ & &
\hspace{1cm}+(1+\varepsilon')\varepsilon +2\varepsilon' C=\\ & &
\langle-L_{xF(x)}v,
L_{\tilde{y_2}F(x)}L_{G(x)\tilde{y}_{2}}L_{G(\tilde{x})G(x)}\circ
dG(\tilde{x})(-L_{x\tilde{x}}v)\rangle\\ &
&\hspace{1cm}+\delta+\varepsilon' C'
+(1+\varepsilon')L+(1+\varepsilon')\varepsilon +2\varepsilon'
C\leq\\ & & \delta+\varepsilon' C'
+(1+\varepsilon')L+(1+\varepsilon')\varepsilon +2\varepsilon' C+
\\ &
&\hspace{1cm}+\langle-L_{xF(x)}v,
L_{\tilde{y_2}F(x)}L_{G(x)\tilde{y}_{2}}L_{G(\tilde{x})G(x)}L_{G(x)G(\tilde{x})}\circ
dG(x)(-v)\rangle +\varepsilon'=
\\ &
&\delta+\varepsilon' C'
+(1+\varepsilon')L+(1+\varepsilon')\varepsilon +2\varepsilon'
C+\varepsilon' +\\ & &\hspace{1cm}+\langle-L_{xF(x)}v,
L_{\tilde{y_2}F(x)}L_{G(x)\tilde{y}_{2}}\circ dG(x)(-v)\rangle \leq
\\ &
&\delta+\varepsilon' C'
+(1+\varepsilon')L+(1+\varepsilon')\varepsilon +2\varepsilon'
C+\varepsilon' +\\ & &\hspace{1cm}+\langle-L_{xF(x)}v,
L_{F(x)G(x)}\circ dG(x)(-v)\rangle +\varepsilon' C\leq
\\ &
&K+\delta+\varepsilon' C'
+(1+\varepsilon')L+(1+\varepsilon')\varepsilon +2\varepsilon'
C+\varepsilon'+\varepsilon' C,
\end{eqnarray*}
that is,
    $$
    \langle\zeta, -v\rangle\leq\mu(\delta,\varepsilon'):=
    K+\delta+\varepsilon' (1+C' )+(1+\varepsilon')L+(1+\varepsilon')\varepsilon
+3\varepsilon' C. \eqno (13)
    $$
Since $\delta$ and $\varepsilon'$ can be chosen to be arbitrarily
small and
    $$
    \lim_{(\delta,\varepsilon')\to (0,0)}\mu(\delta,\varepsilon')=
    K+L+\varepsilon,
    $$
this argument shows that
    $$
    \langle\zeta, -v\rangle\leq K+L+\varepsilon. \eqno (14)
    $$
Finally, this implies that
\begin{eqnarray*}
& &\|v+\zeta\|\geq \langle v, v+\zeta\rangle=\langle v,v\rangle+
\langle v,\zeta\rangle=\\ & &1-\langle \zeta, -v\rangle\geq
1-(K+L+\varepsilon):=\lambda>0,
\end{eqnarray*}
and $\lambda$ is clearly independent of $x,\zeta$. Moreover,
according to assumption $(vi)$, we have that $d(x_{0},
J(x_{0},G(x_{0})))<R(1-(L+K+\varepsilon))=R\lambda$, so we got
all we needed. \qed
%%%%%%%%%%%%%%%%

Finally let us see what Theorem \ref{supertheorem for manifolds}
means when we consider some special cases of the perturbing
mapping $J$. In the general case of a complete Riemannian
manifold, if we have a mapping $G:M\to M$ having an almost fixed
point $x_{0}$ and satisfying certain conditions, and we compose
$G$ with a mapping $H$ which is relatively close to the identity,
we get that $F=H\circ G$ has a fixed point. More precisely, we
have the following.

%%%%%%%%%%%%%%%%%
\begin{thm}\label{corollary for composition}
Let $M$ be a complete Riemannian manifold, and $G:M\to M$ a $C^1$
smooth function such that $G$ is $C$-Lipschitz on a ball
$B(x_{0},R)$. Let $H:M\to M$ be a differentiable mapping. Assume
that $H(G(x))\notin\textrm{sing}(x)\cup\textrm{sing}(G(x))$ for
every $x\in B(x_{0}, R)$, that $$\langle L_{xH(G((x)))}h,
L_{G(x)H(G(x))}dG(x)(h)\rangle_{F(x)}\leq K<1$$ for all $x\in
B(x_{0},R)$ and $h\in TM_{x}$ with $\|h\|_{x}=1$, and that
$\|dH(y)-L_{yH(y)}\|<\varepsilon/C$ for every $y\in G(B(x_{0},
R))$, where $\varepsilon<1-K$, and $d(x_{0}, H(G(x_{0})))<
R(1-K-\varepsilon)$. Then $F=H\circ G$ has a fixed point in
$B(x_{0}, R)$.

If $M$ is finite dimensional one can replace $\textrm{sing}(z)$
with $\textrm{cut}(z)$ everywhere.
\end{thm}
%%%%%%%%%%%%%%%%%
\begin{proof}
It is enough to consider the mapping $J(x,y)=H(y)$. Since
$x\mapsto J_{y}(x)$ is constant for every $y$, we can apply
Theorem \ref{supertheorem for manifolds} with $L=0$.
\end{proof}
%%%%%%%%%%%%%%%%%

Notice that when we take $0<R<\rho=i(M)$, the global injectivity
radius of $M$, we obtain the first Corollary mentioned in the
general introduction.

As another consequence we also have a local version of the
result, whose statement becomes simpler.

%%%%%%%%%%%%%%%%%%
\begin{thm}
Let $M$ be a complete Riemannian manifold. Let $x_0$ be a fixed
point of a $C^1$ function $G:M\to M$ satisfying the following
condition: $$\langle h, dG(x_0)(h)\rangle \leq K<1 \ \ \
\hbox{for every} \ \ \ \|h\|=1.$$ Then there exists a positive
$\delta$ such that for every differentiable mapping $H:M\to M$
such that $\|dH(y)-L_{yH(y)}\|<\delta$
for every  $y$ near $x_0$, the composition $H\circ G:M\to M$ has a fixed
point provided that $d(x_0, H(x_0))<\delta$.
\end{thm}
%%%%%%%%%%%%%%%%%%

If $M$ is endowed with a Lie group structure a natural extension
of Corollary \ref{main fixed point theorem for Hilbert spaces}
holds: we can perturb the function $G$ by summing a small
function $H$ with a small Lipschitz constant, and we get that
$G+H$ has a fixed point.

%%%%%%%%%%%%%%%%%
\begin{thm}\label{fixed point theorem for Lie groups}
Let $(M,+)$ be a complete Riemannian manifold with a Lie group
structure. Let $G:M\to M$ be a $C^1$ smooth function which is
$C$-Lipschitz on a ball $B(x_{0}, R)$. Let $H:M\to M$ be an
L-Lipschitz function. Assume that
$G(x)+H(x)\notin\textrm{sing}(x)\cup\textrm{sing}(G(x))$ for every
$x\in B(x_{0},R)$, and that $$\langle L_{x(H(x)+G(x))}h,
L_{G(x)((H(x)+G(x)))}dG(x)(h)\rangle_{F(x)}\leq K<1$$ for all
$x\in B(x_{0},R)$ and $h\in TM_{x}$ with $\|h\|_{x}=1$.  Then
$G+H$ has a fixed point, provided that $L<1-K$ and $d(x_0,
x_0+H(x_0))<R(1-K-L)$.

Again, if $M$ is finite dimensional one can replace
$\textrm{sing}(z)$ with $\textrm{cut}(z)$ everywhere.
\end{thm}
%%%%%%%%%%%%%%%%%
\begin{proof}
Define $J(x,y)=y+H(x)$. We have that
    $$
    \frac{\partial H}{\partial y}(x,y)=L_{yJ(x,y)},
    $$
so we can apply Theorem \ref{supertheorem for manifolds} with
$\varepsilon=0$.
\end{proof}
%%%%%%%%%%%%%%%%%

Let us conclude with an analogue of Theorem \ref{local fixed point
theorem}, which can be immediately deduced from Theorem \ref{fixed
point theorem for Lie groups}.

%%%%%%%%%%%%%%%%%
\begin{cor}
Let $(M,+)$ be a complete Riemannian manifold with a Lie group
structure. Let $x_0$ be a fixed point of a $C^1$ function $G:M\to
M$ satisfying the following condition: $$\langle h,
dG(x_0)(h)\rangle \leq K<1 \ \ \ \hbox{for every} \ \ \ \|h\|=1.$$
Then there exists a positive $\delta$ such that for every
Lipschitz mapping $H:M\to M$ with Lipschitz constant smaller than
$\delta$, the mapping $G+H:M\to M$ has a fixed point provided that
$d(x_0,x_0+H(x_0))<\delta$.
\end{cor}
%%%%%%%%%%%%%%%%%

%Let us show an easy example of a situation in which the above results are applicable.
%Consider the two-dimensional torus $T^2$ as a square of $\mathbb{R}^2$ centered
%at the origin with opposite edges identified. Let $G:T^2\to T^2$ be the mapping
%induced on $T^2$ by a rotation of $\mathbb{R}^2$ of angle $\pi/2$ and center
%the origin, and let $x_0$ be the image of the origin by the natural quotient
%map from the plane $\mathbb{R}^2$ onto the torus. The mapping $G$ takes the
%Greenwich meridian of $T^2$ onto the exterior equator of $T^2$. It is clear
%that $\langle h, dG(x_0)(h)\rangle=0$ for all $h$, so we can apply the above
%corollaries with $K=0$ to obtain that, if %$H:T^2\to T^2$ is any mapping so
%that $H(x_0)$ is close to $x_0$ and $dH$ is close to $L_{yH(y)}$ %when $y$
%is near $x_0$, then the mapping $H\circ G:T^2\to T^2$ has a fixed point near $x_0$.
%Moreover, if we consider the natural Lie group structure of $T^2$ induced by
%the quotient map from the plane onto $T^2$, then there exists some $\delta>0$
%so that $G+H$ has a fixed point near $x_0$ for %every $\delta$-Lipschitz mapping
%$H\circ G:T^2\to T^2$ so that $H(x_0)$ is $\delta$-close to $x_0$.

Let us show an easy example of a situation in which the above
results are applicable. Let $M$ be the cylinder defined by
$x^2+y^2=1$ in $\mathbb{R}^3$, and let $G:M\to M$ be the mapping
defined by $G(x,y,z)=(x,-y,-z)$. Take $p_{0}$ to be either
$(1,0,0)$ or $(-1,0,0)$ (the only two fixed points of $G$). We
have that $G$ is $1$-Lipschitz and $\langle L_{pq}h,
L_{G(p)q}dG(p)(h)\rangle=-1:=K$ whenever
$q\notin\textrm{cut}(p)\cup \textrm{cut}(G(p))$. Then we can
apply Theorem \ref{corollary for composition} with $R=\pi/2$ to
obtain that, if we take any differentiable mapping $H:M\to M$ such
that $H(G(p))\notin\textrm{cut}(p)\cup\textrm{cut}(G(p))$ for
every $p\in B(p_{0}, \pi/2)$ and
$\|dH(p)-L_{pH(p)}\|<\varepsilon$ for every $p\in G(B(p_{0},R))$
and $d(p_{0}, H(G(p_{0}))<R(1-K-\varepsilon)$, where
$0<\varepsilon<2$, then the composition $F=H\circ G$ has a fixed
point in $B(p_{0},\pi/2)$. 

In a similar way one can also apply Theorem
\ref{fixed point theorem for Lie groups} to obtain that, when $M$ is endowed
with the natural Lie goup structure of $S^{1}\times\mathbb{R}$, the mapping $G+H$
has a fixed point near $p_0$ provided $H:M\to M$ is a relatively small Lipschitz function.

\medskip

\noindent {\bf Acknowledgement}. We thank Francisco R. Ruiz del
Portal for several valuable conversations concerning the results
of this paper.

\medskip

%%%%%%%%%%%%%%%%%%%%%%

\end{document}